\documentclass[12pt]{amsart}
\usepackage[colorlinks]{hyperref}
\hypersetup{
    colorlinks = true,
    linkcolor={blue},
    citecolor={blue},
    urlcolor={black},
}
\usepackage{color,graphicx,shortvrb}
\usepackage[latin 1]{inputenc}

\usepackage[active]{srcltx} 

\usepackage{enumerate}
\usepackage{amssymb}
\usepackage{tikz-cd}
\usepackage{mathrsfs}
\usepackage [ all ]{xy}
\usepackage{soul}
\usepackage{ulem}

\newcommand{\C}{\mbox{${\mathcal C}$}}

\newcommand{\U}{\mbox{${\mathcal U}$}}
\newcommand{\V}{\mbox{${\mathcal V}$}}

\newcommand{\fg}{\mbox{${\mathfrak g}$}}
\newcommand{\fL}{\mbox{${\mathfrak L}$}}

\newcommand{\sspan}{\operatorname{span}}

\newtheorem{theorem}{Theorem}[section]

\newtheorem{corollary}[theorem]{Corollary}

\newtheorem{definition}[theorem]{Definition}
\newtheorem{lemma}[theorem]{Lemma}

\newtheorem{proposition}[theorem]{Proposition}
\newtheorem{remark}[theorem]{Remark}

\newtheorem{example}[theorem]{Example}

\DeclareGraphicsExtensions{.jpg,.pdf,.png,.eps}

\def\sqr#1#2{{\,\vcenter{\vbox{\hrule height.#2pt\hbox{\vrule width.#2pt
height#1pt \kern#1pt\vrule width.#2pt}\hrule height.#2pt}}\,}}
\def\bo{\sqr55\,}
\newcommand{\cl}[1]{\overline{#1}}

\numberwithin{equation}{section}

\begin{document}

\title[Isomorphisms of TKK Lie algebras of JB*-triples]{Isomorphisms of Tits--Kantor--Koecher Lie algebras  of JB*-triples}

\author[M. Cueto-Avellaneda]{Mar{\'i}a Cueto-Avellaneda}
\address[M. Cueto-Avellaneda]{Centro de An{\'a}lise Matem{\'a}tica, Geometria e Sistemas Din\^{a}micos, Departamento de Matem{\'a}tica, Instituto Superior T{\'e}cnico, Universidade de Lisboa, Av. Rovisco Pais 1, 1049-001 Lisbon, Portugal}
\email{emecueto@gmail.com}

\author[L. Oliveira]{Lina Oliveira}
\address[Lina Oliveira]{Centro de An{\'a}lise Matem{\'a}tica, Geometria e Sistemas Din\^{a}micos, Departamento de Matem{\'a}tica, Instituto Superior T{\'e}cnico, Universidade de Lisboa, Av. Rovisco Pais 1, 1049-001 Lisbon, Portugal
}
\email{lina.oliveira@tecnico.ulisboa.pt}

\thanks{The first author was partially supported by grant PID2021-122126NB-C31 funded by MICIU/AEI/10.13039/501100011033 and by ERDF/EU, and by FCT/Portugal through project UIDB/04459/2020 with DOI identifier 10-54499/UIDP/04459/2020.
The second author was partially supported by FCT/Portugal through project UIDB/04459/2020 with DOI identifier 10-54499/UIDP/04459/2020.}

\subjclass[2020]{ 17B70, 17C50, 17C65, 18A05, 46L70}
\keywords{Tits--Kantor--Koecher Lie algebra, JB*-triple,  isomorphism, isometry, tripotent, equivalence of categories}

\date{\today}

\begin{abstract} We characterise the isomorphisms of Tits--Kantor--Koecher Lie algebras of JB*-triples as a class of surjective linear isometries and  show how these algebras form a category  equivalent to that of JB*-triples. We introduce the concepts of tripotent,  and orthogonality and order amongst tripotents for Tits--Kantor--Koecher Lie algebras. This leads to showing that
a graded or negatively graded order isomorphism between certain subsets of tripotents of two  Tits--Kantor--Koecher Lie algebras of atomic JB*-triples, which commutes with involutions, preserves orthogonality and is continuous at a non-zero tripotent of a specific type,  can be extended as a real-linear isomorphism between the algebras.
\end{abstract}

\maketitle
\thispagestyle{empty}

\section{Introduction}
The existence of a one-to-one correspondence between Jordan triples and a class of Lie algebras, known as the {\it Tits--Kantor--Koecher Lie algebras}, is well-documented in the literature (\cite{k1,k2, Koecher, my, tits}),  with this correspondence being established through the   Tits--Kantor--Koecher construction. This correspondence has attracted significant  attention over the years (see, for example, \cite{CavSmir,Gar,Jacob,loos1,Neher}) and, more recently, the Tits--Kantor--Koecher Lie algebras corresponding to JB*-triples, a class of Jordan triples,  have been identified (\cite{ChuOliveira}).
 In fact, with a careful choice of morphisms, it follows that the category of canonical non-degenerate Tits-Kantor-Koecher Lie algebras is equivalent to that of non-degenerate Jordan triples (\cite[Theorem 1.3.11]{ChuBook}). In light of the identification in \cite{ChuOliveira}, one is led to consider what this equivalence might imply for JB*-triples and their corresponding Tits--Kantor--Koecher Lie algebras. 
 
 This  is precisely the stance we adopt in the present work: we explore the strong  categorical connection between JB*-triples and their Tits--Kantor--Koecher Lie algebras (TKK Lie algebras, for short). 
  In Section 2, we discuss in detail the equivalence between the category of JB*-triples and that of their TKK Lie algebras. Building on this equivalence, we  prove a central result (Theorem \ref{t isomorphisms}) which  characterises the isomorphisms of Tits--Kantor--Koecher Lie algebras of JB*-triples as a class of surjective linear isometries.

Having characterised isomorphisms in Section 2, we approach them from a different perspective in Section 3. In that section, we introduce the concepts of tripotent,  and orthogonality and order amongst tripotents for  Tits--Kantor--Koecher Lie algebras and investigate their properties. This leads to the proof of Theorem \ref{t_atomic}, the main result of this section, which shows  that
 a graded  order isomorphism  between certain subsets  of  tripotents of two  Tits--Kantor--Koecher Lie algebras of atomic JB*-triples, which commutes with involutions, preserves  orthogonality and is continuous at a non-zero  tripotent of a specific type,  can be extended as a real-linear  isomorphism between the algebras. For negatively graded mappings, a corresponding statement is provided in Corollary \ref{c_atomic}.

We finish Section 1 by recalling some facts needed in the following sections.\\

A {\it Hermitian Jordan triple} $V$ is a complex vector space equipped with a {\it triple product}, that is, a mapping  
 $$\{\cdot,\cdot,\cdot\} \colon V \times V \times V \rightarrow V$$
 $$ (a,b,c)\mapsto \{a,b,c\}
 $$
 which is symmetric and linear in the outer variables, conjugate linear in the middle variable  and satisfies, for all $a,b,x,y,z\in V$, the  {\it Jordan triple identity} 
 \begin{equation}\label{eq Jordan triple identity}
     \{a,b,\{x,y,z\}\} = \{\{a,b,x\}, y,z\} - \{x, \{b,a,y\},z\} + \{x,y, \{a,b,z\}\}.
 \end{equation}

A Hermitian Jordan triple $V$ is called {\it non-degenerate} when, for each $a\in V$, if $\{a,x,a\} =0$ for all $x\in V$, then $a=0$. Given two elements $a,b$ in a Hermitian Jordan triple $V$, define the linear operator  $a\bo b : V \rightarrow V$, said a {\it box operator}, by  
$a \bo b(x) = \{a,b,x\},$ 
for all $x\in V$. 

In \cite{Ka83}, W. Kaup introduced a  class of Hermitian Jordan triples, that of  JB*-triples. A Hermitian Jordan triple $V$ is said to be a {\it JB*-triple}  if $V$ admits additionally the structure of a complex Banach space for which the triple product is continuous and such that, for all $a\in V$,
\begin{enumerate}[{\rm (i)}]
\item $a\bo a$ is a Hermitian operator with non-negative
	spectrum,
\item $\|a\bo a\| = \|a\|^2$.
\end{enumerate}

Examples of JB*-triples are C*-algebras,  JB*-algebras, and the  Cartan factors, the latter classified in six types.
Let $H$ and $K$ be complex Hilbert spaces. The {\it Cartan factors of type I, II and III} are, respectively,  $B(H,K)$, $A(H)=\{z\in B(H): z^t=-z\}$, and $S(H)=\{z\in B(H): z^t=z\}$, with $z^t=Jz^*J$ denoting the {\it transpose} of $z$ in the JB*-triple $B(H)$, where $J\colon H\to H$ is a conjugation, that is, $J$ is a conjugate linear isometric involution.

The remaining types of Cartan factors are the {\it spin factors} (type IV) and the {\it exceptional Cartan factors} $M_{1,2}(\mathbb{O})$ (type V), and  $H_3(\mathbb{O})$ (type VI) of $1\times 2$ matrices and  $3\times 3$ Hermitian matrices over the octonions, respectively. A {\it spin factor} V is a closed subspace of $B(H)$ for some Hilbert
space $H$, such that $v\in V$ implies $v^*\in V$ and $v^2\in \mathbb{C}I_H$.

An element $e$ in a JB*-triple $V$ is said to be a {\it tripotent} if $\{e,e,e\} =e$. For each tripotent $e\in V$, there exists an algebraic decomposition of $V$ known as the {\it Peirce decomposition}\label{eq Peirce decomposition} associated with $e$ determined by the eigenspaces of the operator $e\bo e$. Namely, 
$$V= V_{2} (e) \oplus V_{1} (e) \oplus V_0 (e),
$$ where $V_i (e)=\{ x\in V : \{e,e,x\} = \frac{i}{2}x \}$,
for each $i=0,1,2$. It is easy to see that every Peirce subspace $V_i(e)$ is a JB*-subtriple of $V$.

This splitting of $V$ into the direct sum of the Peirce spaces is subjected to the so-called {\it Peirce arithmetic}:
 $\{V_{i}(e),V_{j} (e),V_{k} (e)\}\subset V_{i-j+k} (e)$ if $i-j+k \in \{ 0,1,2\},$ and $\{V_{i}(e),V_{j} (e),V_{k} (e)\}=\{0\}$ otherwise, and $$\{V_{2} (e),V_{0}(e),V\} =  \{V_{0} (e),V_{2}(e),V\} =\{0\}.$$
The projection $P_{i_{}}(e)$ of $V$ onto $V_{i} (e)$,  $i=0,1,2$, is called the {\it Peirce $i$-projection}.  Peirce projections are contractive (\cite[Corollary 1.2]{FriRu85}) and satisfy the equalities  
$$P_{2}(e) = Q(e)^2,\quad P_{1}(e) =2(e\bo e-Q(e)^2), \quad P_{0}(e) =\hbox{Id}_V - 2 e\bo e + Q(e)^2,$$
where $Q(e):V\to V$ is the conjugate linear map defined, for all $x\in V$, by $Q(e) (x) =\{e,x,e\}$.

Given two elements $a,b$ in a JB*-triple $V$, we  say that $a$ and $b$ are {\it orthogonal} (and we write $a\perp b$), if $a\bo b =0$. The following characterisations hold:
$$a\perp b \Leftrightarrow \{a,a,b\} =0 \Leftrightarrow \{b,b,a\}=0 \Leftrightarrow b\perp a.$$ 

Let $e$ be a tripotent in $V$. It can be deduced from Peirce arithmetic that any two elements $a$ and $b$ in $V$ are orthogonal whenever $a\in V_2(e)$ and $b\in V_0(e)$.
For tripotents $e,u$ in $V$, we have that $$ e\perp u \Leftrightarrow e\pm u \hbox{ are tripotents}$$ (see \cite[Lemma 3.6]{IsKaRo95}).

Let $\U(V)$ be the set of all tripotents in a JB*-triple $V$. For
two elements $e$ and $u$ of $\mathcal U(V)$, write $u \leq e$, 
if $e - u$ is a tripotent orthogonal to $u$. 
This defines a partial ordering on $\mathcal U(V)$ (see \cite{batt, ERFS}).

The existence of non-zero tripotents  is not guaranteed in a JB*-triple. Notwithstanding, there is a class of JB*-triples possessing an abundance of tripotents, that of  JBW*-triples. 
 A {\it JBW*-triple} is a JB*-triple which is also a dual Banach space with a unique isometric predual (\cite{BarTi86}). It is well-known that the second dual of a JB*-triple is a JBW*-triple (\cite{Di86}), and that  the triple product of a JBW*-triple is separately weak*-continuous  (\cite{BarTi86, Horn87}). 

It is also the case that the extreme points of the closed unit ball of a JBW*-triple consist entirely of tripotents. Consequently, the 
Krein--Milman theorem guarantees that any  JBW*-triple is well-supplied with tripotents (see \cite[Lemma 4.1]{BraKaUp78}, \cite[Theorem 3.2.3]{ChuBook}, \cite[Proposition 3.5]{KaUp77}).
The weak*-closed linear span of the family $\mathcal U(V)$ of tripotents in a JBW*-triple $V$ is in fact  equal to $V$ itself (\cite{FriRu85}). 
Moreover, the set $\widetilde{\mathcal U}(V)=\U\cup \{\omega\}$ of tripotents in a JBW*-triple $V$  with a greatest element adjoined is a complete lattice, when endowed with the partial ordering above  (see \cite{batt, ERFS}).

An important property of JB*-triples, which we will rely frequently on in the sequel, is that when two such objects are isomorphic they are isometrically so.  
Let $V,W$ be JB*-triples. A mapping $\varphi\colon V\to W$ is said to be a {\it triple isomorphism}  if $\varphi$ is a linear bijection preserving the triple product, that is, for all $a,b,c\in V$,
$$\varphi (\{a,b,c\})=\{\varphi(a), \varphi(b), \varphi(c)\}.$$
Kaup's Banach--Stone theorem states that a linear bijection between JB*-triples is an isometry if and only if it is a triple isomorphism (see \cite[Proposition 5.5]{Ka83}).

It was shown in \cite{ChuOliveira} that there exists a correspondence between JB*-triples and a specific class of Tits--Kantor--Koecher Lie algebras. We  describe now this correspondence for future reference. To do so, we  recall firstly the definition of  those Lie algebras.

Let $\mathfrak g$ be a real or complex Lie algebra and let  $[\cdot, \cdot]:\mathfrak{g}\times \mathfrak{g}\to \mathfrak{g}$ be its Lie multiplication. The Lie algebra  $\mathfrak g$ is said to be $3$-{\it graded} if it is  a direct sum 
$$\mathfrak g = \mathfrak g_{-1} \oplus \mathfrak g_0 \oplus \mathfrak g_1$$ 
of linear subspaces 
satisfying $[\mathfrak g_m, \mathfrak g_n] \subset \mathfrak g_{m+n}$, where $\mathfrak g_{m+n} = \{0\}$, whenever $m+n \notin \{-1,0, 1\}$. 
Furthermore, if $[\mathfrak g_{-1}, \mathfrak g_1]=\mathfrak g_0$, then $\mathfrak g$ is said to be {\it canonical}.

An {\it involution} on a Lie algebra $\mathfrak{g}$ is an automorphism $\theta\colon  \mathfrak g \to \mathfrak g$ such that $\theta^2=\mbox{I}_\mathfrak{g}$. Here, it is understood that $\theta$ is conjugate linear, if $\mathfrak{g}$ is a complex Lie algebra. A Lie algebra $(\mathfrak{g},\theta)$ with an involution is said to be {\it involutive}.

\begin{definition}
An involutive Lie algebra  $(\mathfrak g, \theta)$ is said to be  a {\rm Tits--Kantor--Koecher Lie algebra}  if $\mathfrak g =\mathfrak g_{-1}\oplus\mathfrak g_0\oplus \mathfrak g_1$ is a $3$-graded Lie algebra and its involution  $\theta\colon \mathfrak g\to \mathfrak g$ is {\it negatively graded}, that is,  $\theta(\mathfrak g_j)=\mathfrak g_{-j}$, for $j = 0, \pm 1$. 
\end{definition}
For the sake of brevity, henceforth a Tits--Kantor--Koecher Lie algebra might be referred to as a {\it TKK Lie algebra}, and the corresponding involution $\theta$ will be called the TKK involution of $\mathfrak{g}$. 
In what follows, we consider only TKK involutions and might refer to them simply as involutions.
We occasionally denote the TKK Lie algebra $(\mathfrak g, \theta)$  by $\mathfrak g$, for convenience.

Let $V$ be a complex vector space, we  use the symbol $\cl{V}$ to denote the {\it conjugate} of $V$, that is, a complex vector space given by the set $V$ equipped with its original addition operation but a new scalar multiplication defined, for all $\lambda\in \mathbb{C}, x\in \cl{V}$, by
\begin{align*}
    \cdot :\mathbb{C}\times \cl{V}&\longrightarrow \cl{V}\\
    (\lambda,x)&\longmapsto \lambda\cdot x=\cl{\lambda}x.
\end{align*}

Observe that, if $V$ is a JB*-triple, then  $\cl{V}$ is also a JB*-triple with the same triple product and norm (see \cite{Dang92, Ka83}).

We describe next the Tits--Kantor--Koecher construction of the TKK Lie algebra $\mathfrak{L}(V)$ associated with a   Hermitian Jordan triple $V$. 
Given a Hermitian Jordan triple $V$, consider the direct sum 
\begin{equation}\label{eq constrution}
    \mathfrak L (V):=V\oplus V_0\oplus\cl{V},
\end{equation}
where $\cl{V}$ is the conjugate of $V$ and $V_0= \sspan\{a\bo b: a,b \in V\}$ denotes the linear span of the box operators on $V$.  For  $(x,h,y)$ and $(u,k,v)$ in $\fL(V)=V\oplus V_0\oplus\cl{V}$, define
\begin{equation}\label{eq Lie bracket}
[(x, h, y), (u,k,v)]=(h(u)-k(x),  [h,k]+x\bo v-u\bo y, k^\natural(y)-h^\natural(v)),
\end{equation}
where $\natural: V_0 \rightarrow V_0$ is an involutive conjugate linear mapping  defined on elements  $h = \sum_j a_j \bo b_j \in V_0$ by $h^\natural = \sum_j b_j \bo a_j$.

The space $\mathfrak{L}(V)$ together with the bracket \eqref{eq Lie bracket} is  a $3$-graded Lie algebra $\mathfrak L (V) = \mathfrak L (V)_{-1} \oplus \mathfrak L (V)_0 \oplus \mathfrak L (V)_1$, where
$$\mathfrak L (V)_{-1} = V, \quad \mathfrak L (V)_0 = V_0\quad \text{and } \mathfrak L (V)_1 = \cl{V}.
$$
Moreover, $\mathfrak{L}(V)$ is an involutive Lie algebra for the negatively graded involution $\theta\colon V\oplus V_0\oplus\cl{V}\to V\oplus V_0\oplus\cl{V}$ defined, for all  $(x,h,y) \in  V\oplus V_0\oplus\cl{V}$, by
$$\theta(x, h,y)=(y,-h^\natural, x).
$$
It is clear that $(\mathfrak L (V), \theta)$ is a canonical Tits-Kantor-Koecher Lie algebra. We shall frequently refer to $\fL(V)$ as the TKK Lie algebra of $V$ or associated with $V$.

It is worth pointing out that, for $a,b,c\in V$, one has 
\begin{equation}\label{eq triple product in Lie language}
    \{a,b,c\}= [[a,\theta(b)],c].
\end{equation}

Notice that, for the  TKK Lie algebra $(\mathfrak{L}(V), \theta)$ of a Hermitian Jordan triple $V$, one has 		
\begin{equation}
    a\bo b = \left[ a, \theta b \right]=-\theta(b\bo a), \quad \forall a,b\in \mathfrak{L}_{-1}(V).
\end{equation}
Furthermore, if $V$ is a JB*-triple, then its  TKK Lie algebra $(\mathfrak{L}(V),\theta)$  is a {\it normed Lie algebra} for the norm 
\begin{equation}\label{eq norm}
\|(x,h,y)\| = \|x\| + \|h\| +\|y\|, \qquad \left((x,h,y) \in V \oplus V_0\oplus \cl{V}\right).\end{equation}
In other words, 
$\mathfrak{L}(V)$ is a normed vector space and
  $$  \|[X,Y]\| \leq C \|X\|\|Y\| \qquad (X,Y \in \mathfrak{L}(V)),
  $$
for some constant $C>0$.  The continuity of the Jordan triple product on $V$ implies that the TKK involution $\theta$ is continuous. 

Let $ad$ be the adjoint representation of a TKK Lie algebra $\mathfrak g= \mathfrak g_{-1}\oplus\mathfrak g_0\oplus \mathfrak g_1$.
For $a \in \mathfrak g_{-1}$ and $b\in \mathfrak g_{1}$, we have $[a,b] \in \mathfrak g_0$ and hence the adjoint $ad\,[a,b] : \mathfrak g \rightarrow \mathfrak g$  is a {\it graded} derivation, that is, $ad\,[a,b](\mathfrak g_\ell) \subset \mathfrak g_\ell$
for $\ell =0, \pm 1$. A TKK Lie algebra $\mathfrak g$ is said to be {\it non-degenerate} if, for $a\in \mathfrak g_{1}$,
$$(ad\, a)^2=0 \implies a=0.$$
The TKK Lie algebra $\fL(V)$ of a non-degenerate Hermitian Jordan triple is a non-degenerate TKK Lie algebra.

The next theorem characterises the TKK Lie algebras of  JB*-triples.

\begin{theorem}\cite[Theorem 3.4]{ChuOliveira}\label{t ChuOliveira}
 Let $\mathfrak g$ be a complex Lie algebra. The following conditions are equivalent.
 \begin{itemize}
\item[(i)] $\mathfrak g$ is the TKK Lie algebra of a JB*-triple.
\item[(ii)] $\mathfrak g$ has a real form $\mathfrak g_r$, which is isomorphic to the reduced Lie algebra of a bounded symmetric domain.
\item[(iii)] $\mathfrak g$ is a normed canonical $3$-graded Lie algebra $\mathfrak g_{-1}\oplus \mathfrak g_0\oplus\mathfrak g_1$ with a negatively graded continuous involution $\theta$ such that $\mathfrak g_{-1}$ is a Banach space in the given norm and each $a \in \mathfrak g_{-1}$ with $\|a\|=1$ satisfies
 \begin{itemize}
\item[(a)] $\|[[a, \theta a], a]\| = 1$,
\item[(b)] $i {\rm ad}\, [a, \theta a](z) \frac{\partial}{\partial z} \in  \mathfrak h$, where $\mathfrak h$ is the Lie algebra of $$\{x \in \mathfrak g_{-1}: \|x\|<1\},$$
\item[(c)] $I- {\rm ad}\, [a, \theta a] : \mathfrak g_{-1} \rightarrow \mathfrak g_{-1}$ is not invertible, where $I$ is the identity operator on $\mathfrak g$. 
\end{itemize}
\end{itemize}
\end{theorem}
\noindent
Here, the real form  $\mathfrak g_r$ is a real subalgebra of the Lie algebra $\mathfrak g_{\mathbb{R}}$, which is  $\mathfrak g$ considered as a Lie algebra over the reals (cf. \cite{Helgason}, p.180), such that 
$\fg=\fg_r+i\fg_r$.\\

In Section 3, we prove a theorem about the extension of certain mappings on atomic JBW*-triples, that is, JBW*-triples that are realised as $\ell_\infty$-sums of Cartan factors (cf. Theorem \ref{t_atomic}). More precisely, a JBW*-triple $V$ is said to be {\it atomic} if $V= \bigoplus_{\alpha\in \Lambda}^\infty  V_\alpha$,  where $\Lambda$ is some set of indices and all summands $V_\alpha$ are Cartan factors.

For completeness, we recall briefly the definition of the $\ell_\infty$-sum of a family of JB*-triples and some facts about its TKK Lie algebra (cf., \cite{ChuOliveira}).

Let   $\{V_\alpha\}_{\alpha\in \Lambda}$ be a family of JB*-triples, and let 
$$V=\bigoplus_{\alpha \in \Lambda}^\infty V_\alpha:= \{ (v_\alpha) \in \Pi_\alpha V_\alpha
: \sup_\alpha \|v_\alpha\| < \infty\}
$$ 
be their $\ell_\infty$-sum.
Recall that $V$ is a JB*-triple where the triple product is defined  coordinatewise and the $\ell_\infty$-norm is defined by 
$$\|(v_\alpha)\| : =\sup_\alpha \|v_\alpha\|.
$$
We can define the direct sum of the box operators $a_\alpha \bo b_\alpha: V_\alpha \rightarrow V_\alpha$,
for $(a_\alpha), (b_\alpha) \in \oplus_\alpha ^\infty V_\alpha$, by
$$\oplus_\alpha (a_\alpha \bo b_\alpha): (v_\alpha) \in \bigoplus_\alpha^\infty  V_\alpha \mapsto
( a_\alpha \bo b_\alpha (v_\alpha)) \in \bigoplus_\alpha^\infty  V_\alpha$$
and we have
$$(a_\alpha) \bo (b_\alpha) = \oplus_\alpha (a_\alpha \bo b_\alpha).$$
Consequently, the linear span  $(\bigoplus_\alpha^\infty  V_\alpha)_0$ of the box operators on
$\bigoplus_\alpha^\infty  V_\alpha$ is contained in the $\ell_\infty$-sum $\bigoplus_\alpha^\infty ( V_\alpha)_0$,
where $(V_\alpha)_0$ is the linear span of box operators on $V_\alpha$.  
Moreover,
$$\mathfrak L(V) = \bigoplus_\alpha^\infty  V_\alpha \oplus (\bigoplus_\alpha^\infty  V_\alpha)_0 \oplus
\bigoplus_\alpha^\infty  \cl{V_\alpha}\subset \bigoplus_\alpha^\infty  V_\alpha \oplus \bigoplus_\alpha^\infty  (V_\alpha)_0 \oplus
\bigoplus_\alpha^\infty  \cl{V_\alpha} = \mathfrak L,$$
where 
$\mathfrak L:=\bigoplus_\alpha^\infty \mathfrak L (V_\alpha)$ is also a $3$-graded Lie algebra with grading
$$\mathfrak L_{-1} = \oplus_\alpha^\infty V_\alpha, \quad  \mathfrak L_{0} = \oplus_\alpha^\infty (V_\alpha)_0,
\quad \mathfrak L_{1} = \oplus_\alpha^\infty \cl{V_\alpha}.$$

\begin{lemma}\cite[Lemma 2.1]{ChuOliveira}\label{dsum} Let $V= \bigoplus_\alpha^\infty  V_\alpha$ be the $\ell_\infty$-sum of a family
$\{V_\alpha\}_{\alpha \in \Lambda}$ of JB*-triples. Then the TKK Lie algebra $\mathfrak L(V)$ of $V$ is contained in the
$\ell_\infty$-sum $\bigoplus_\alpha^\infty \mathfrak L (V_\alpha)$ of the TKK Lie algebras $(\mathfrak L(V_{\alpha}))_\alpha$ and the TKK involution of $\mathfrak L(V)$
is (the restriction of) the direct
sum of the TKK involutions of $(\mathfrak L(V_{\alpha}))_\alpha$.
\end{lemma}

\section{Functoriality}\label{s Funct}

Non-degenerate TKK Lie algebras form a category whose morphisms are the graded isomorphisms (that is, graded linear bijections preserving the Lie product) commuting with involutions. This category is known to be equivalent to that consisting of non-degenerate Jordan triples and triple isomorphisms (\cite[Theorem 1.3.11]{ChuBook}). 

Let   $\V$ be the category of JB*-triples whose objects are JB*-triples and morphisms are triple isomorphisms. This is a subcategory of the category of non-degenerate Jordan triples.
On the other hand,  the subset of TKK Lie algebras associated with the JB*-triples via the TKK construction (see \eqref{eq constrution} and Theorem \ref{t ChuOliveira}) forms a category $\C$ where the objects are those TKK Lie algebras characterised in Theorem \ref{t ChuOliveira} and the morphisms are the graded isomorphisms commuting with involutions. This category $\C$ will be called the {\it category of TKK Lie algebras of JB*-triples}. It is also the case that $\C$ is a subcategory of the category of non-degenerate TKK Lie algebras.  

In Theorem \ref{t isomorphisms}, we characterise the isomorphisms between TKK Lie algebras of JB*-triples as a specific class of surjective isometries, and our approach to its proof relies on the equivalence of the categories $\V$ and $\C$. For this reason and to make this section self-contained, before proving Theorem \ref{t isomorphisms}, we will describe in detail the corresponding equivalence functor. For the general theory of categories, see, for example, \cite{Adamek, MacLane, Riehl}.  

Define the mapping $F:\mathcal{V}\to \mathcal{C}$, respectively, on objects $V\in \V$ and morphisms $\varphi\colon V\to W$ ($V,W\in \V$),  by  
\begin{equation}\label{02}
V\mapsto FV:=\fL(V)
\end{equation}
and 
\begin{equation}
\label{03}F\varphi: FV\to FW
\end{equation}
\begin{equation}\label{04}
    F\varphi(x,h,y):= \left(\varphi(x),\varphi\circ h\circ \varphi^{-1}, \varphi(y)\right), \qquad\quad  (x,h,y)\in FV.
\end{equation}

\begin{lemma}\label{l_functor}
Let $\V$ be the category of JB*-triples and let $\C$ be the category of TKK Lie algebras of JB*-triples. Then, the mapping $F:\mathcal{V}\to \mathcal{C}$ defined in \eqref{02}--\eqref{04} is a functor.
\end{lemma}

\begin{proof}
It follows from the TKK construction and Theorem \ref{t ChuOliveira} that $FV\in \mathcal{C}$. Recall that we also have the grading  
$$FV=\mathfrak{L}(V)= \mathfrak{L}(V)_{-1}\oplus \mathfrak{L}(V)_0\oplus \mathfrak{L}(V)_1= V\oplus V_0\oplus \cl{V}.$$

We begin by showing that the mapping in \eqref{03},\eqref{04} is well-defined. Given $x\in FV_{-1}=V$, we have that $\varphi(x)\in W= FW_{-1}$. On the other hand, for $y\in FV_1=\cl{V}$, we have that $y$ lies also in the $V$, consequently, $\varphi(y)\in\cl{W}=FW_{1}$, as $W$and $\cl{W}=FW_{1}$ coincide as sets.

Finally, let $h=\sum_{j=1}^{n} a_j\bo b_j$ lie in $FV_0=V_0$, where $a_j,b_j\in V$, for all $j=1,\dots, n$.  It follows that, for all $w\in W$,
\begin{align*}
    \varphi \circ h \circ \varphi^{-1} (w)=& \varphi\bigl(\sum_j a_j\bo b_j(\varphi^{-1}(w))\bigr)=\varphi(\sum_j\{ a_j,b_j,\varphi^{-1}(w) \})\\ =& \sum_j \varphi(\{a_j, b_j, \varphi^{-1}(w)\}) = \sum_j \{\varphi(a_j), \varphi(b_j), w \} \\ =& \sum_j \varphi(a_j)\bo \varphi(b_j)(w).
\end{align*}
Hence, 
\begin{equation}\label{05}
\varphi \circ h \circ \varphi^{-1}=\varphi(\sum_{j=1}^{n} a_j\bo b_j)\varphi^{-1}
= \sum_j \varphi(a_j)\bo \varphi(b_j),
\end{equation}
yielding that $\varphi \circ h \circ \varphi^{-1}$ is a finite sum of box operators on $W$  which, therefore, lies in $W_0$.
 We may now conclude that $F\varphi$ is graded. It is also easy to see from the definition that $F\varphi$ is a linear bijection. We show now that $F\varphi$ commutes with the involutions.

Let $\theta, \theta'$ be the involutions on $V$ and $W$, respectively,  and let $(x,h,y)\in FV$. Then, by \eqref{05},
$$
F\varphi\theta(x,h,y)=F\varphi(y,-h^\natural,x)
= (\varphi(y),-(\varphi(h))^\natural,\varphi(x))=\theta'F\varphi.
$$

To see that $F\varphi$ is a graded isomorphism commuting with involutions, it only remains to show that $F\varphi$ preserves the Lie bracket.
For  $(x,h,y),(u,k,v)\in FV$, the definition  \eqref{eq Lie bracket} of the Lie product leads to 
$$
F\varphi[(x, h, y), (u,k,v)]=F\varphi\left( h(u)-k(x), [h,k]+x\bo v-u\bo y, k^\natural(y)-h^\natural(v)\right).
$$
Hence, 
$$
F\varphi[(x, h, y), (u,k,v)]= (\varphi (h(u)-k(x)), \varphi ([h,k]+x\bo v-u\bo y) \varphi^{-1}, \varphi(k^\natural(y)-h^\natural(v))).
$$
Setting $h=\sum_{j=1}^{n}a_j\bo b_j$ and $k=\sum_{l=1}^{m}c_l\bo d_l$ with  $a_j,b_j, c_l, d_l\in V$, for all $j=1,\dots, n$ and  $l=1,\dots, m$, we have
\begin{align*}
\varphi (h(u)-k(x))=&\varphi(\sum_{j}a_j\bo b_j(u) - \sum_{l}c_l\bo d_l(x) )\\=&\sum_{j}\varphi(a_j)\bo \varphi(b_j)(\varphi(u)) - \sum_{l}\varphi(c_l)\bo \varphi(d_l)(\varphi(x) )\\=&
\varphi \circ h \circ \varphi^{-1}(\varphi(u))-\varphi\circ k\circ \varphi^{-1}(\varphi(x)).
\end{align*}

Similarly,
\begin{align*}
    \varphi(k^\natural(y)-h^\natural(v))&= \varphi (\sum_{l=1}^{m}d_l\bo c_l(y)-\sum_{j=1}^{n}b_j\bo a_j(v))\\
    &=\sum_{l=1}^{m}\varphi (d_l)\bo \varphi (c_l)(\varphi (y))-\sum_{j=1}^{n}\varphi (b_j)\bo \varphi (a_j) (\varphi (v))\\
    &=(\varphi\circ k\circ \varphi^{-1})^\natural(\varphi(y))-(\varphi\circ h\circ \varphi^{-1})^\natural(\varphi(v)).
\end{align*}

As to the middle term, we have 
\begin{align*}
\varphi\circ ([h,k]+x\bo v-u\bo y)\circ \varphi^{-1}=& \varphi\circ [h,k]\circ \varphi^{-1}+\varphi\circ x\bo v\circ \varphi^{-1}\\&-\varphi\circ u\bo y\circ \varphi^{-1}\\=& \varphi hk \varphi^{-1} - \varphi kh\varphi^{-1}+\varphi (x)\bo \varphi(v)\\&-\varphi(u)\bo \varphi(y)\\=& [\varphi\circ h \circ \varphi^{-1},\varphi\circ k\circ \varphi^{-1}] + \varphi (x)\bo \varphi(v)\\&-\varphi(u)\bo \varphi(y).
\end{align*}
We have just shown that $F\varphi$ preserves the Lie bracket, that is, for all $(x,h,y),(u,k,v)\in FV$,
$$F\varphi [(x,h,y),(u,k,v)] = [F\varphi(x,h,y), F\varphi(u,k,v)].
$$ 

All of the above leads to the conclusion that  $F\colon \V\to \C$ maps arrows of $\V$ to arrows of $\C$.
Moreover, let $\varphi:V\to W$ and $\varphi':W\to V'$ be two triple isomorphisms between the JB*-triples $V,W$ and $V'$, and  consider the composition mapping $\varphi'\circ \varphi: V\to V'$. By  \eqref{03}, \eqref{04}, it is easily seen  that $F(\varphi'\circ \varphi)=F\varphi'\circ F\varphi$. It is also the case that, for each $V\in \mathcal{V}$, we have $F\hbox{Id}_V=\hbox{Id}_{FV}$, where $\hbox{Id}_{V}$ and $\hbox{Id}_{FV}$ denote the identity maps on $V$ and $FV$, respectively.

We have just shown that $F:\mathcal{V}\to \mathcal{C}$ is a functor between the categories $\mathcal{V}$ and $\mathcal{C}$.
\end{proof}

\begin{theorem}\label{t_equivalence}
    Let $\mathcal{V}$ be the category of JB*-triples, and let   $\mathcal{C}$ be the category of TKK Lie algebras of JB*-triples. Then,  $F:\mathcal{V}\to \mathcal{C}$ defined in \eqref{02}--\eqref{04}  is an equivalence of categories.
\end{theorem}

\begin{proof} By Lemma \ref{l_functor}, we know that $F:\mathcal{V}\to \mathcal{C}$ is a functor. It suffices to show that the functor $F$ is full, faithful, and essentially surjective on objects to guarantee that $F$ is an equivalence of the categories $\mathcal{V}$ and $\mathcal{C}$ (see \cite[Definition 3.33]{Adamek},\cite[Theorem 1, p.91]{MacLane}, \cite[Theorem 1.5.9]{Riehl}). 
     
We begin by showing that $F$ is full, that is, given a morphism $T:FV\to FW$ in the category $\mathcal{C}$, for some $V,W\in \V$, there exists a triple isomorphism $\varphi:V\to W$  such that $T=F\varphi$. 
     
Define $\varphi:V\to W$ by $\varphi= T_{|_{V}}$, where $V$ is identified with $V\oplus\{0\}\oplus\{0\}$. The mapping $\varphi$ is a triple isomorphism, since $T$ is a graded isomorphism commuting with involutions. Then,
$$
 F\varphi(x,h,y)= (\varphi(x), \varphi\circ h\circ \varphi^{-1}, \varphi(y))= (T(x), \varphi\circ h\circ \varphi^{-1}, T(y) ).
$$
Let  $\theta,\theta'$ be the involutions on $FV,FW$, respectively, and  let $h=\sum_{j=1}^{n} a_j\bo b_j$, with $a_j,b_j\in V$, for all $j=1,\dots, n$. Since $T$ is a morphism in $\C$, we have
\begin{align*}
        T(h)=& T(\sum_j a_j\bo b_j)=\sum_j T(a_j\bo b_j)\\=& \sum_jT[a_j,\theta b_j]= \sum_j [T(a_j), \theta'T(b_j)]\\=& \sum_j T(a_j)\bo T(b_j)=\sum_j \varphi(a_j)\bo \varphi (b_j) = \varphi\circ h\circ \varphi^{-1}.
\end{align*}
Hence, $F$ is full. Furthermore, suppose there exist triple isomorphisms $\varphi$ and $\varphi'$  such that $F\varphi = T = F\varphi'$. It follows immediately that  $\varphi=\varphi'$, yielding that $F$ is faithful. 

By \cite[Theorem 1.3.11]{ChuBook} and Theorem \ref{t ChuOliveira}, there exists a one-to-one correspondence between JB*-triples and the objects of the category $\C$. Hence, it follows immediately that $F$ is essentially surjective on objects. We have established that $F$ is an equivalence of categories.
\end{proof}

The following corollary is an immediate consequence of this theorem and of Theorem 5.4 in \cite{Ka83}.

\begin{corollary}
The category of all bounded symmetric domains with base point is equivalent to the category of TKK Lie algebras of JB*-triples.
\end{corollary}

The next theorem establishes a relation between isometries of TKK Lie algebras of JB*-triples and the arrows of the category $\C$.

\begin{theorem}\label{t isomorphisms}
Let $V, W$ be  JB*-triples, and let $(\mathfrak{L}(V),\theta),(\mathfrak{L}(W), \theta')$ be their TKK Lie algebras, respectively. Let $T: \mathfrak{L}(V)\to \mathfrak{L}(W)$ be a graded linear mapping commuting with involutions. Then, $T$ is an isomorphism if and only if $T$ is a surjective isometry. 
\end{theorem}
\begin{proof}
Suppose that $T: \mathfrak{L}(V)\to \mathfrak{L}(W)$ is an isomorphism.  Then, by Theorem \ref{t_equivalence},  
there exists a triple isomorphism $\varphi\colon V\to W$ such that $T=F\varphi$. Hence, by Kaup's Banach--Stone Theorem (\cite[Proposition 5.5]{Ka83}), $T_{|_V}$ is a surjective isometry, and similarly for $T_{|_{\overline{V}}}$. Recall that, by definition, $V=\mathfrak{L}(V)_{-1}$ and $\overline{V}=\mathfrak{L}(V)_{1}$.

Now, let $h=\sum_{j=1}^n a_j\bo b_j$, with $a_j,b_j\in \mathfrak{L}(V)_{-1}$, for all $j=1,\dots, n$. Then, 
    \begin{align*}
        ||T(h)||&=||T(\sum_{j=1}^n a_j\bo b_j)||=||T(\sum_{j=1}^n [a_j,\theta b_j])||\\&=||\sum_{j=1}^n [T(a_j),\theta' T(b_j)]||=||\sum_{j=1}^n T(a_j)\bo T(b_j)||.
    \end{align*}
Keeping in mind  that $T$ is a bijection and a triple isomorphism on $\mathfrak{L}(V)_{-1}$, we have 
    \begin{align*}
        ||T(h)||&= \sup_{v\in\mathfrak{L}(W)_{-1}, \|v\|=1} ||\sum_{j=1}^n T(a_j)\bo T(b_j)(v)||\\
        &=  \sup_{y\in\mathfrak{L}(V)_{-1}, \|y\|=1} ||\sum_{j=1}^n T(a_j)\bo T(b_j)(T(y))||\\
        &=  \sup_{y\in\mathfrak{L}(V)_{-1}, \|y\|=1} ||T(\sum_{j=1}^n \{a_j, b_j,y\})||.
    \end{align*}
Since $T|_{\mathfrak{L}(V)_{-1}}$ is a linear isometry, it now follows that  
    \begin{align*}||T(h)||&=  \sup_{y\in\mathfrak{L}(V)_{-1}, \|y\|=1} ||T(\sum_{j=1}^n \{a_j, b_j,y\})||\\
      &=  \sup_{y\in\mathfrak{L}(V)_{-1}, \|y\|=1} ||\sum_{j=1}^n \{a_j, b_j,y\}||\\
       &=  \sup_{y\in\mathfrak{L}(V)_{-1}, \|y\|=1} ||\sum_{j=1}^n a_j\bo b_j(y)||
       =||\sum_{j=1}^n a_j\bo b_j||=||h||.
    \end{align*}
Hence, for all $(x,h,y)\in\fL(V)$,
   $$
    \|T(x,h,y)\|=\|T(x)\|+\|T(h)\|+\|T(y)\|
    =\|x\|+\|h\|+\|y\|=\|(x,h,y)\|,
   $$
which shows that $T$ is an isometry, as required.\smallskip

Suppose now that $T: \mathfrak{L}(V)\to \mathfrak{L}(W)$ is a surjective isometry. Then, $T_{|_V}$ determines a triple isomorphism $\varphi\colon V\to W$. Hence, by Theorem \ref{t_equivalence}, there exists a (iso)morphism 
    $F\varphi\colon \mathfrak{L}(V)\to \mathfrak{L}(W)$. We will show that $F\varphi$ and $T$ coincide. 

    Let $T_{j}:=T|_{\mathfrak{L}(V)_{j}}$ and $(F\varphi)_{j}:=(F\varphi)|_{\mathfrak{L}(V)_{j}}$ for $j=0,\pm 1$. Notice that, by the definition of $F$, it is clear that 
    $(F\varphi)_{-1}=T_{-1}$ (see \eqref{02}-\eqref{04}). On the other hand, since $T$ commutes with involutions, we have, for $y\in \mathfrak{L}(V)_{1}$,
    $$\theta'T_{-1}(\theta y)=(\theta')^2T_{1}( y)=T_{1}( y).
    $$
     Hence, $T_{1}=\theta'(F\varphi)_{-1}=(F\varphi)_{1}.$

Let $a,b\in \mathfrak{L}(V)_{-1}=V$ and consider the operator $a\bo b$. 
 For $x\in \mathfrak{L}(V)_{-1}$, we have that there exists uniquely $y\in \mathfrak{L}(W)_{-1}$ such that
$$T(a\bo b)(x)=T\{a,b,x\}=(F\varphi)_{-1}\{a,b,x\}=\varphi\{a,b,\varphi^{-1}(y)\}.
$$
Hence, 
$T(a\bo b)=\varphi\circ a\bo b\circ \varphi^{-1}.$
It follows from linearity that also, for $h\in\mathfrak{L}(V)_0$, 
$T(h)=\varphi\circ h\circ \varphi^{-1}.$
 Observing that $\mathfrak{L}(V)$ is a canonical TKK Lie algebra, we have finally that $T_0=(F\varphi)_0$, concluding the proof.
\end{proof}

\begin{theorem}\label{t isomorphisms conj}
    Let $V, W$ be  JB*-triples, and let $(\mathfrak{L}(V),\theta),(\mathfrak{L}(W), \theta')$ be their TKK Lie algebras, respectively. Let $T: \mathfrak{L}(V)\to \mathfrak{L}(W)$ be a graded conjugate-linear mapping commuting with involutions. 
    If $T$ is a conjugate-linear isomorphism then $T$ is a surjective isometry. 
\end{theorem}
\begin{proof}
    Suppose $T$ is a (graded) conjugate-linear isomorphism. Since $\mathfrak{L}(V)_{-1}=V$ and $\mathfrak{L}(W)_{-1}=W$, equality \eqref{eq triple product in Lie language} can be applied to deduce that the conjugate-linear bijection given by the restriction $T|_{\mathfrak{L}(V)_{-1}}:\mathfrak{L}(V)_{-1}\to \mathfrak{L}(W)_{-1}$ preserves the triple product. Indeed, for any $a,b,c\in V$,
    \begin{align*}
    T|_V(\{a,b,c\}_V)=&T([[a,\theta b], c])= [[T(a), \theta ' T(b)], T(c)]\\=&\{T(a),T(b),T(c)\}_W.    
    \end{align*}
    It follows that $T|_{\mathfrak{L}(V)_{-1}}$ is a conjugate-linear surjective isometry. 
    
    Similarly, it can be seen that the conjugate-linear bijection $T|_{\mathfrak{L}(V)_{1}}:\mathfrak{L}(V)_{1}\to \mathfrak{L}(W)_{1}$ preserves the triple product on $\mathfrak{L}(V)_{-1}=\cl{V}$: for any $\cl{a}, \cl{b}, \cl{c}\in \cl{V}$,

   $$ T|_{\cl{V}}(\{\cl{a},\cl{b},\cl{c}\}_{\cl{V}})=  T_{|_{\cl{V}}}(\{\theta_{|_{V}}(a),\theta_{|_{V}}(b),\theta_{|_{V}}(c)\}_{\cl{V}})\\
    = T_{|_{\cl{V}}}\circ \theta_{|_{V}}(\{a, b, c\}_{V}).$$
    By \eqref{eq triple product in Lie language} and noticing that $T$ commutes with involutions, it follows that
    \begin{align*}
    T_{|_{\cl{V}}}(\{\cl{a},\cl{b},\cl{c}\}_{\cl{V}})= & T\circ \theta([[a,\theta b],c]) \\
    =& \theta ' \circ T ([[a,\theta b],c])\\
    =& [[T(\theta(a)), \theta'T(\theta(b))], T(\theta (c))]\\
    =& \{ T_{|_{\cl{V}}}(\cl{a}), T_{|_{\cl{V}}}(\cl{b}), T_{|_{\cl{V}}}(\cl{c}) \}_{\cl{W}}.
    \end{align*}
    Therefore, we have that $T|_{\mathfrak{L}(V)_{1}}:\mathfrak{L}(V)_{1}\to \mathfrak{L}(W)_{1}$ is a conjugate-linear surjective isometry.\smallskip

    The same arguments used in the proof of Theorem \ref{t isomorphisms} can be applied here to prove that $T|_{\mathfrak{L}(V)_{0}}:\mathfrak{L}(V)_{0}\to \mathfrak{L}(W)_{0}$ is a surjective isometry.\smallskip

    Collecting all the partial results obtained above, we have that, for all $(x,h,y)\in \mathfrak{L}(V)$,  
    \begin{align*}
        || T(x,h,y) ||&= || T|_{\mathfrak{L}(V)_{-1}}(x) || + ||T|_{\mathfrak{L}(V)_{0}}(h) ||+ || T|_{\mathfrak{L}(V)_{1}}(y) || \\&= || x|| +|| h||+||y ||=||(x,h,y) ||.
    \end{align*}
        
\end{proof}

\section{Tripotents}
Tripotents in JB*-triples play a crucial role in the theory of these spaces (see, for instance, \cite{CuetoBecFerPer, EdFerHosPe2010, EdRu96, ERFS,Sidd2007}). Bearing in mind the equivalence of categories binding JB*-triples and a class of Tits--Kantor--Koecher Lie algebras (Theorem \ref{t_equivalence}), one is naturally led to ponder the existence of corresponding elements in the latter. The aim of the present section is precisely  that: 
we  introduce the concept of tripotent in a Tits--Kantor--Koecher Lie algebra and explore the related notions of orthogonality and order amongst these special elements. We  subsequently study the linear extension of a particular class of  bijections defined on a subset of tripotents of the TKK Lie algebra of a JB$^*$-triple.
In order to facilitate writing the calculations, we now set some notation. 

Given a TKK Lie algebra $(\mathfrak{g}, \theta)$, for each $j\in\{-1,0,1\}$, we  write $P^\mathfrak{g}_{j}$ to denote each  of the natural projections from $\mathfrak{g}$ to $\mathfrak{g}_{-1}\oplus\{0\}\oplus\{0\}$,  $\{0\}\oplus\mathfrak{g}_{0}\oplus\{0\}$, and $\{0\}\oplus\{0\}\oplus\mathfrak{g}_{1}$, respectively. That is, given $z=(x,h,y)\in \mathfrak{g}$, we write $P^\mathfrak{g}_{-1}(z):=(x,0,0)$, $P^\mathfrak{g}_{0}(z):=(0,h,0)$, and $P^\mathfrak{g}_{1}(z):=(0,0,y)$. This notation allows for writing $z\in \mathfrak{g}$ as 
$$z=P^\mathfrak{g}_{-1}(z)+ P^\mathfrak{g}_{0}(z) + P^\mathfrak{g}_{1}(z).
$$
Clearly, given $z\in \mathfrak{g}$, we have, for all $j=-1,0,1$, that $$P^\mathfrak{g}_j(\theta z)=\theta P^\mathfrak{g}_{-j}(z).
$$
It will also be useful for future calculations to note that, for any two elements $z,w\in \mathfrak{g}$, $[P^\mathfrak{g}_{j}(z), P^\mathfrak{g}_{j}(w)]=0$ for $j=-1,1$.

Let $\mathfrak{g}_{\pm1}$ be the subspace of $\mathfrak{g}$ defined by
$\mathfrak{g}_{\pm1}=\{z\in\mathfrak{g}\colon P^\mathfrak{g}_0(z)=0\}. 
$

\begin{definition}\label{d tripotent}
Let $\mathfrak{g}=\mathfrak{g}_{-1}\oplus\mathfrak{g}_0\oplus \mathfrak{g}_1$ be a TKK Lie algebra with involution $\theta$. An element $z\in \mathfrak{g}$ is said to be a {\rm tripotent} if 
$$\left[ \left[  z, \theta z  \right] , z \right]  = z.$$  
The set of all tripotents in $\mathfrak{g}$ is denoted by $\mathcal{U}(\mathfrak{g})$. The notation $\mathcal{U}(\mathfrak{l})$ will refer to the set of all the elements of a subset $\mathfrak{l}\subseteq \mathfrak{g}$ which are tripotents in $\mathfrak{g}$.
\end{definition}

\begin{example}
A JB*-triple $V$ is {\it flat} if, for all $a,b\in V$,   $a\bo b=b\bo a$. Given a flat JB*-triple $V$, the element $(0,a\bo b,0) \in \fL(V)$ cannot be a tripotent unless it is zero.

In fact, this also happens in abelian JB*-triples. A JB*-triple $V$ is said to be  {\it abelian} if, for all $a,b,c,d\in V$, one has $[a\bo b,c\bo d]=0$. A direct application of Definition \ref{d tripotent} shows that  $(0,a\bo b, 0)\in \fL(V)$ is a tripotent if and only if $a\bo b=0$ (that is, $a\perp b$ in the JB$^*$-triple $V$).

Clearly, for a (not necessarily flat or abelian) JB*-triple $V$, the element  $(0,a\bo a,0)\in \fL(V)$ is a tripotent  if only if  $a=0$.
\end{example}

\begin{example}

Let $A$ be a unital C$^*$-algebra and, for $a,b\in A$, consider $z=(0, a\bo b, 0)\in \mathfrak{L}(A)$. Then, 
$$
[z, \theta z]= [(0,a\bo b, 0),(0, -b\bo a, 0)] = (0,[a\bo b, -b\bo a],0),
$$
 where, for any  $c\in A$,
\begin{align*}
[a\bo b, -b\bo a](c)&= -[a\bo b, b\bo a](c) \\
&= -(a\bo b)\circ (b\bo a) (c)+ (b\bo a)\circ (a\bo b) (c) \\
&= \frac14 (  ba^*ab^*c + ba^*cb^*a+ ab^*ca^*b + cb^*aa^*b \\&- ab^*ba^*c - ab^*ca^*b - ba^*cb^*a - ca^*bb^*a   ).
\end{align*}

If $a$ and $b$ are unitary elements in $A$, that is, if  $aa^*=1_A=a^*a$ and $b^*b=1_A=bb^*$, then
\begin{align*}
[a\bo b, -b\bo a](c)&= \frac14 (c + ba^*cb^*a+ ab^*ca^*b + c\\&- c - ab^*ca^*b - ba^*cb^*a - c   )=0.
\end{align*}
It follows from Definition \ref{d tripotent} that $z$ is a tripotent in $\mathfrak{L}(A)$ if and only if $a\bo b=0$ (or, in other words,  $a\perp b$ in $A$).
\end{example}

In view of the examples above , we will give a particular attention to tripotents $z\in \mathcal{U}(\mathfrak{g})$ in a TKK Lie algebra $\mathfrak{g}$ such that $P^\mathfrak{g}_0(z)=0$. Specifically, we consider  the set $\U(\mathfrak{g_{\pm1}})$   of elements in $\mathfrak{g}_{\pm1}$ which are tripotents in $\mathfrak{g}$, and the set $\U_s(\mathfrak{g_{\pm1}})$ of elements $z\in \U(\mathfrak{g_{\pm1}})$ such that, for $j=-1,1$,  
\begin{equation}\label{eq strict Lie tripotents}
		\left[ \left[  P^\mathfrak{g}_j(z), \theta P^\mathfrak{g}_j(z)  \right] , P^\mathfrak{g}_j(z) \right]  = P^\mathfrak{g}_j(z).
\end{equation}

According to the above, we have the following chain of inclusions for a TKK Lie algebra  $\mathfrak{g}$:
$$
\U_s(\mathfrak{g_{\pm1}})\subseteq\U(\mathfrak{g_{\pm1}})\subseteq \U(\mathfrak{g}).
$$

Notice that the involution in a TKK Lie algebra $(\mathfrak{g}, \theta)$ preserves tripotents. Indeed, given $z\in \mathcal{U}(\mathfrak{g})$, since the involution $\theta$ preserves the Lie brackets, we have that 
$$
\left[ \left[  \theta z, \theta(\theta z)  \right] , \theta z \right]  = \theta \left[ \left[  z, \theta z  \right] , z \right]  = \theta z,
$$
yielding that $\theta z\in \mathcal{U}(\mathfrak{g})$.  Since $\mathfrak{g}_{\pm1}$ is invariant under $\theta$, it immediately follows that the involution also preserves elements in $\mathcal{U}(\mathfrak{g}_{\pm 1})$. 
Moreover, if $z\in \mathcal{U}_s(\mathfrak{g}_{\pm 1})$, we have, for each $j=-1,1$, 
\begin{align*}
    \left[ \left[  P^\mathfrak{g}_j(\theta z), \theta P^\mathfrak{g}_j(\theta z)  \right] , P^\mathfrak{g}_j(\theta z) \right] =& \left[ \left[  \theta P^\mathfrak{g}_{-j}(z), P^\mathfrak{g}_{-j}(z)  \right] , \theta P^\mathfrak{g}_{-j}(z) \right]\\
    =& \theta \left[ \left[  P^\mathfrak{g}_{-j}(z), \theta P^\mathfrak{g}_{-j}(z)  \right] , P^\mathfrak{g}_{-j}(z) \right] \\
    =& \theta P^\mathfrak{g}_{-j}(z) = P^\mathfrak{g}_{j}(\theta z),
\end{align*}
showing that the set   $\mathcal{U}_s(\mathfrak{g}_{\pm 1})$ is also left invariant under the involution $\theta$.

We collect some properties derived from the Lie-algebraic notion of general tripotents and from those in $\mathcal{U}_s(\mathfrak{g}_{\pm1})$. We will frequently use the fact that identities \eqref{eq strict Lie tripotents} might be satisfied by general elements in $\fg$.

\begin{proposition}\label{l properties of strictness}
    Let $(\mathfrak{g},\theta)$ be a TKK Lie algebra and let $z\in \mathfrak{g}_{\pm1}$. Suppose $z$ satisfies identities \eqref{eq strict Lie tripotents}. Then,
    \begin{itemize}
        \item[(i)] $P^\mathfrak{g}_j(z)\in \mathcal{U}_s(\mathfrak{g_{\pm1}})$, for $j=-1,1$;
        \item [(ii)]There exist $z_1,z_2\in \mathcal{U}_s(\mathfrak{g}_{\pm1})$ such that $z=z_1+z_2$;
        \item[(iii)] If $z\in \mathcal{U}_s(\mathfrak{g_{\pm1}})$, then there exist $z_1,z_2\in\mathcal{U}_s(\mathfrak{g}_{\pm1})$ such that,  for $j=-1,1$, $z=P^\mathfrak{g}_j(z_1)+\theta P^\mathfrak{g}_j(z_2)$.
\end{itemize}
\end{proposition}

\begin{proof}

    (i)  For $j=-1,1$, it is clear that $P^\mathfrak{g}_j(z)\in \mathcal{U}(\mathfrak{g}_{\pm1})$ by \eqref{eq strict Lie tripotents}. Moreover, for $k=-1,1$, we have that
    \begin{align*}
        \left[ \left[  P^\mathfrak{g}_k(P^\mathfrak{g}_j(z)), \theta P^\mathfrak{g}_k(P^\mathfrak{g}_j(z))  \right] , P^\mathfrak{g}_k(P^\mathfrak{g}_j(z)) \right]  = 0 = P^\mathfrak{g}_k(P^\mathfrak{g}_j(z)),
    \end{align*}
    whenever $k\neq j$, and by \eqref{eq strict Lie tripotents},
    \begin{align*}
        \left[ \left[  P^\mathfrak{g}_k(P^\mathfrak{g}_j(z)), \theta P^\mathfrak{g}_k(P^\mathfrak{g}_j(z))  \right] , P^\mathfrak{g}_k(P^\mathfrak{g}_j(z)) \right]  &= \left[ \left[  P^\mathfrak{g}_j(z), \theta P^\mathfrak{g}_j(z)  \right] , P^\mathfrak{g}_j(z) \right]\\  &= P^\mathfrak{g}_j(z) = P^\mathfrak{g}_k(P^\mathfrak{g}_j(z)),
    \end{align*}
    provided that $k=j$. We have just shown that $P^\mathfrak{g}_j(z)\in \mathcal{U}_s(\mathfrak{g}_{\pm1})$ for $j=-1,1$.\smallskip

     (ii) By (i) of this lemma, it is enough to consider $z_1=P^\mathfrak{g}_{-1}(z)$ and $z_2=P^\mathfrak{g}_1(z)$.\smallskip

     (iii) It is enough to consider $z_1=z$ and $z_2=\theta z$.
\end{proof}

An observation that will be useful later is that the definition of tripotent given in Definition \ref{d tripotent}  for TKK Lie algebras aligns with the notion of a tripotent for JB*-triples, when working with their associated TKK Lie algebras.

\begin{proposition}\label{p tripotent characterisation}
Let $V$ be a JB*-triple, let $\mathfrak{L}(V)=\mathfrak{L}(V)_{-1}\oplus\mathfrak{L}(V)_0\oplus \mathfrak{L}(V)_1$ be its  TKK Lie algebra and let $z$ be an element in $V$. 
Then, $z$ is a tripotent in $V$ if and only if $z$ is a tripotent in $\mathfrak{L}(V)$. 
\end{proposition}
\begin{proof}
By the TKK construction, we know that  $z$ is identified with $P^{\mathfrak{L}(V)}_{-1}(z)\in \mathfrak{L}(V)$, and $P^{\mathfrak{L}(V)}_{0}(z)=0=P^{\mathfrak{L}(V)}_{1}(z)$. That is to say that $z\in \mathfrak{g}_{\pm1}$.
        
Now, suppose firstly that $z$ is a tripotent in $V$. Applying \eqref{eq triple product in Lie language}, it is straightforward to see that $z\in \mathcal{U}(\mathfrak{L}(V))$. 
On the other hand, if $z$ is a  tripotent in $\mathfrak{L}(V)$, then, by    \eqref{eq triple product in Lie language}, 
$$
z=\left[ \left[  z, \theta z  \right] , z \right]=\{z,z,z\},
$$
which shows that $z$ is a tripotent in the JB*-triple $V$.
\end{proof}

Observe that, in the conditions of Proposition \ref{p tripotent characterisation}, any tripotent $z\in V$ of a JB$^*$-triple $V$ will automatically be an element of the set $\mathcal{U}_s(\mathfrak{L}(V)_{\pm1})$, when seen as an element in the associated TKK Lie algebra $\mathfrak{L}(V)$.

The characterisation of the tripotents of a JB*-triple as tripotents in the corresponding TKK Lie algebra naturally leads to the study of whether some of the triple properties of these elements could be translated to a Lie-algebraic form. For instance, we are interested in defining a partial order on the set of tripotents. In order to do that, we begin with the concept of orthogonality.

\begin{definition}\label{d orth}
Let $\mathfrak{g}=\mathfrak{g}_{-1}\oplus\mathfrak{g}_0\oplus\mathfrak{g}_1$ be a TKK Lie algebra with involution $\theta$.
Given  $z,w$  in $\mathfrak{g}_{\pm 1}$,  $z$ is said to be {\rm orthogonal to} $w$, written $z\perp w$,  if
$[ P^\mathfrak{g}_j(z), \theta P^\mathfrak{g}_j(w) ] =0,$ for $j=-1,1$.
\end{definition}

The notion of orthogonality in the triple setting coincides with that introduced in the Lie algebra context when working with TKK Lie algebras associated with JB*-triples, as shown in the following proposition.

\begin{proposition}\label{p orth characterisation}
Let $V$ be a JB$^*$-triple, let $\mathfrak{L}(V)=\mathfrak{L}(V)_{-1}\oplus\mathfrak{L}(V)_0\oplus \mathfrak{L}(V)_1$ be its TKK Lie algebra with involution $\theta$, 
 and let $z,w\in V$.Then,  $z$ and $w$ are orthogonal in $V$ if and only if $z$ and $w$ are orthogonal in $\mathfrak{L}(V)$. 
\end{proposition}

\begin{proof}
    Given $z,w\in V$, since $V=\mathfrak{L}(V)_{-1}$, we have that $z$ and $w$ can be identified with $P^\mathfrak{g}_{-1}(z)$ and $P^\mathfrak{g}_{-1}(w)$, respectively. 
    
    Firstly, suppose that $z$ and $w$ are (triple) orthogonal in $V$, that is, $z\bo w=0$. Hence, by \eqref{eq Lie bracket},
    $$[ P^\mathfrak{g}_{-1}(z), \theta P^\mathfrak{g}_{-1}(w) ] = z\bo w=0.
    $$ 
    On the other hand, $P^\mathfrak{g}_{1}(z)=0= P^\mathfrak{g}_{1}(w)$ from which follows that $z$ and $w$ are orthogonal.

    To prove the converse, suppose that $z$ and $w$ are orthogonal in $\mathfrak{L}(V)$. Then, by definition, $[ P^\mathfrak{g}_j(z), \theta P^\mathfrak{g}_j(w) ] =0,$  for $j=-1,1$. In particular, 
    $$0=[ P^\mathfrak{g}_{-1}(z), \theta P^\mathfrak{g}_{-1}(w) ] =z\bo w,
    $$ yielding  that $z$ and $w$ are orthogonal in the JB*-triple $V$.    
\end{proof}

Let $(\mathfrak{g}, \theta)$ and $(\mathfrak{h}, \theta')$ be two TKK Lie algebras, and let $\varphi: \mathfrak{g} \to \mathfrak{h}$ be a mapping such that $\varphi(\mathfrak{g}{\pm 1}) \subseteq \mathfrak{h}{\pm 1}$. Note that a specific example of such a mapping is any graded, or negatively graded, $\varphi: \mathfrak{g} \to \mathfrak{h}$ that sends zero to zero in $\mathfrak{g}_0$. We say that $\varphi$ preserves orthogonality (in one direction) if, for any two orthogonal elements $z, w \in \mathfrak{g}_{\pm 1}$, we have $\varphi(z) \perp \varphi(w)$ in $\mathfrak{h}$. That is, if $[ P^\mathfrak{g}_j(z), \theta P^\mathfrak{g}_j(w) ] = 0$, for any $j = -1, 1$, then it follows that $[ P^\mathfrak{h}_j(\varphi(z)), \theta' P^\mathfrak{h}_j(\varphi(w)) ] = 0$ for any $j = -1, 1$.

\smallskip

Let $(\mathfrak{g}, \theta)$ be a TKK Lie algebra. It is worth observing that the TKK Lie involution $\theta$ preserves orthogonality. Indeed, given two elements $z,w\in \mathfrak{g}_{\pm 1}$ with $z\perp w$,  we have, for  $j=-1,1$, that
$$[ P^\mathfrak{g}_j(\theta z), \theta P^\mathfrak{g}_j(\theta w) ] = \theta [ P^\mathfrak{g}_{-j}(z), \theta P^\mathfrak{g}_{-j}(w) ]=\theta(0)=0.
$$
Therefore, $\theta z\perp \theta w.$ 

\begin{remark}\label{remx2}
The TKK involution of a TKK Lie algebra is not the only example of a mapping preserving orthogonality between TKK Lie algebras. In fact, any  
graded homomorphism $T:(\mathfrak{g}, \theta)\to (\mathfrak{h}, \theta ')$ between two TKK Lie algebras $(\mathfrak{g}, \theta)$ and $(\mathfrak{h}, \theta')$ preserves orthogonality provided that $T$ commutes with the involutions. Namely, for any two  orthogonal elements $z,w\in \mathfrak{g}_{\pm 1}$, we have that 
    \begin{align*}
        [ P^\mathfrak{h}_j(T(z)), \theta' P^\mathfrak{h}_j(T(w)) ] &=  [ T(P^\mathfrak{g}_{j}(z)), \theta' T(P^\mathfrak{g}_{j}(w)) ]= [ T(P^\mathfrak{g}_{j}(z)), T\theta (P^\mathfrak{g}_{j}(w)) ]\\&=T[ P^\mathfrak{g}_{j}(z), \theta P^\mathfrak{g}_{j}(w) ]=T(0)=0,
    \end{align*}
    for  $j=-1,1$. Hence, $T(z)\perp T(w)$. 
    
    Note that we would reach the same conclusion by considering $T:(\mathfrak{g}, \theta)\to (\mathfrak{h}, \theta ')$ as a negatively graded homomorphism that commutes with involutions.
\end{remark}

The following lemma is a simple consequence of the Lie-algebraic definition of orthogonality.
\begin{lemma}\label{l about orth basics}
Let $\mathfrak{g}=\mathfrak{g}_{-1}\oplus\mathfrak{g}_0\oplus \mathfrak{g}_1$ be a TKK Lie algebra with involution $\theta$. The following assertions hold.
\begin{itemize}
    \item[(i)] For any $z,w\in \mathfrak{g}$, $P^\mathfrak{g}_j(z)\perp P^\mathfrak{g}_{-j}(w)$, for $j=-1,1$.
    \item[(ii)] For any $z,w\in \mathfrak{g}_{\pm1}$, $z\perp w$ if and only if $w\perp z$.
    \item[(iii)] For any $z,w\in \mathfrak{g}_{\pm1}$, $z\perp w$ if and only if $P^\mathfrak{g}_j(z)\perp P^\mathfrak{g}_j(w)$ for $ j=-1,1.$ 
\end{itemize}
\qed
\end{lemma}

\begin{lemma}\label{l perp}
Let $\mathfrak{g}=\mathfrak{g}_{-1}\oplus\mathfrak{g}_0\oplus \mathfrak{g}_1$ be a TKK Lie algebra with involution $\theta$, and let $z,w\in \mathfrak{g}_{\pm 1}$ be such that $z\perp w$. The following holds.
\begin{itemize}
\item[(i)] For   $j=-1,1$, 
$$[[ P^\mathfrak{g}_j(z), \theta P^\mathfrak{g}_j(z) ], P^\mathfrak{g}_j(w) ]  =0=[[ P^\mathfrak{g}_j(w), \theta P^\mathfrak{g}_j(w) ], P^\mathfrak{g}_j(z) ].$$

\item[(ii)] If $z,w\in \mathcal{U}_s(\mathfrak{g}_{\pm 1})$, 
then $P^\mathfrak{g}_{j}(z)+ \theta P^\mathfrak{g}_{j}(w)\in \mathcal{U}_s(\mathfrak{g}_{\pm 1})$, for $j=-1,1$. 

\end{itemize}

\end{lemma}

\begin{proof}
    (i) By the Jacobi identity, for $j=-1,1$, we have
    \begin{align*}
        [[ P^\mathfrak{g}_j(z), \theta P^\mathfrak{g}_j(z) ], P^\mathfrak{g}_j(w) ]=& - [[ \theta P^\mathfrak{g}_j(z), P^\mathfrak{g}_j(w) ], P^\mathfrak{g}_j(z) ]\\&-[[ P^\mathfrak{g}_j(w), P^\mathfrak{g}_j(z) ], \theta P^\mathfrak{g}_j(z) ] \\=&
        [[ P^\mathfrak{g}_j(w), \theta P^\mathfrak{g}_j(z) ], P^\mathfrak{g}_j(z) ] = 0.
    \end{align*}
The remaining assertion is proved similarly.

(ii) Since  $z,w\in \mathcal{U}_s(\mathfrak{g}_{\pm 1})$, it is clear that, for each $j=-1,1$, the element $P^\mathfrak{g}_{j}(z)+ \theta P^\mathfrak{g}_{j}(w)$ lies in $\mathfrak{g}_{\pm 1}$, and also that
    \begin{align*}
         [P^\mathfrak{g}_{j}(z)+ \theta P^\mathfrak{g}_{j}(w), \theta P^\mathfrak{g}_{j}(z)+ P^\mathfrak{g}_{j}(w)]&
         = [P^\mathfrak{g}_{j}(z), \theta P^\mathfrak{g}_{j}(z)]+ [\theta P^\mathfrak{g}_{j}(w), P^\mathfrak{g}_{j}(w)].
    \end{align*}
Hence, by orthogonality, it follows that
     \begin{align*}
         &\left[[P^\mathfrak{g}_{j}(z)+ \theta P^\mathfrak{g}_{j}(w), \theta P^\mathfrak{g}_{j}(z)+ P^\mathfrak{g}_{j}(w)],P^\mathfrak{g}_{j}(z)+ \theta P^\mathfrak{g}_{j}(w)\right] \\
         =& \left[[P^\mathfrak{g}_{j}(z), \theta P^\mathfrak{g}_{j}(z)],P^\mathfrak{g}_{j}(z)+ \theta P^\mathfrak{g}_{j}(w)\right]+ \left[[\theta P^\mathfrak{g}_{j}(w), P^\mathfrak{g}_{j}(w)],P^\mathfrak{g}_{j}(z)+ \theta P^\mathfrak{g}_{j}(w)\right]\\
         =&
         \left[[P^\mathfrak{g}_{j}(z), \theta P^\mathfrak{g}_{j}(z)],P^\mathfrak{g}_{j}(z)\right]+ \left[[P^\mathfrak{g}_{j}(z), \theta P^\mathfrak{g}_{j}(z)], \theta P^\mathfrak{g}_{j}(w)\right]\\
         &+ \left[[\theta P^\mathfrak{g}_{j}(w), P^\mathfrak{g}_{j}(w)],P^\mathfrak{g}_{j}(z)\right] + \left[[\theta P^\mathfrak{g}_{j}(w), P^\mathfrak{g}_{j}(w)], \theta P^\mathfrak{g}_{j}(w)\right]\\
         =& P^\mathfrak{g}_{j}(z) + \theta P^\mathfrak{g}_{j}(w) + \left[[P^\mathfrak{g}_{j}(z), \theta P^\mathfrak{g}_{j}(z)], \theta P^\mathfrak{g}_{j}(w)\right] + \left[[\theta P^\mathfrak{g}_{j}(w), P^\mathfrak{g}_{j}(w)],P^\mathfrak{g}_{j}(z)\right]\\
         =& P^\mathfrak{g}_{j}(z) + \theta P^\mathfrak{g}_{j}(w) -\left[[\theta P^\mathfrak{g}_{j}(w), P^\mathfrak{g}_{j}(z)], \theta P^\mathfrak{g}_{j}(z)\right] - \left[[P^\mathfrak{g}_{j}(z), \theta P^\mathfrak{g}_{j}(w)],P^\mathfrak{g}_{j}(w)\right]\\
         =&   P^\mathfrak{g}_{j}(z) + \theta P^\mathfrak{g}_{j}(w),
    \end{align*}
   showing that $P^\mathfrak{g}_{j}(z)+ \theta P^\mathfrak{g}_{j}(w)\in \mathcal{U}(\mathfrak{g})$, for  $j=-1,1$.

    Finally, we consider the element $P^\mathfrak{g}_{k}(P^\mathfrak{g}_{j}(z)+ \theta P^\mathfrak{g}_{j}(w))$ for $k=-1,1$. Suppose $k=j$, then since $P^\mathfrak{g}_{j}(P^\mathfrak{g}_{j}(z)+ \theta P^\mathfrak{g}_{j}(w))=P^\mathfrak{g}_{j}(z)$ and $z\in \mathcal{U}_s(\mathfrak{g}_{\pm 1})$, it follows that, for each $j=-1,1$,
    \begin{align*}
        &[[P^\mathfrak{g}_{j}(P^\mathfrak{g}_{j}(z)+ \theta P^\mathfrak{g}_{j}(w)), \theta P^\mathfrak{g}_{j}(P^\mathfrak{g}_{j}(z)+ \theta P^\mathfrak{g}_{j}(w))], P^\mathfrak{g}_{j}(P^\mathfrak{g}_{j}(z)+ \theta P^\mathfrak{g}_{j}(w))] \\
        &= [[P^\mathfrak{g}_{j}(z), \theta P^\mathfrak{g}_{j}(z)], P^\mathfrak{g}_{j}(z)] = P^\mathfrak{g}_{j}(z)=  P^\mathfrak{g}_{j}(P^\mathfrak{g}_{j}(z)+ \theta P^\mathfrak{g}_{j}(w)).
    \end{align*}

    On the other hand, if $k=-j$, we can conclude that
    \begin{align*}
        &[[P^\mathfrak{g}_{j}(P^\mathfrak{g}_{j}(z)+ \theta P^\mathfrak{g}_{j}(w)), \theta P^\mathfrak{g}_{j}(P^\mathfrak{g}_{j}(z)+ \theta P^\mathfrak{g}_{j}(w))], P^\mathfrak{g}_{j}(P^\mathfrak{g}_{j}(z)+ \theta P^\mathfrak{g}_{j}(w))] \\
        &= [[P^\mathfrak{g}_{j}(z), \theta P^\mathfrak{g}_{j}(z)], P^\mathfrak{g}_{j}(z)] = P^\mathfrak{g}_{j}(z)=  P^\mathfrak{g}_{j}(P^\mathfrak{g}_{j}(z)+ \theta P^\mathfrak{g}_{j}(w)),
    \end{align*}

\end{proof}

Some interesting identities for TKK Lie algebras can be derived from an iterative application of the Jacobi identity.
\begin{lemma}\label{p Lie triple identity} 
Let $(\mathfrak{g},\theta)$ be a TKK Lie algebra. For $j=-1,0,1$, the following identities hold:
\begin{itemize}
    \item[(i)] For any $a,b,x,y,z\in \mathfrak{g}$,
    \begin{align*}
        &\left[\left[ P^\mathfrak{g}_j(a), \theta P^\mathfrak{g}_j(b) \right], \left[\left[ P^\mathfrak{g}_j(x), \theta P^\mathfrak{g}_j(y) \right], P^\mathfrak{g}_j(z)\right]\right] \\
         &= \left[\left[\left[\left[ P^\mathfrak{g}_j(a), \theta P^\mathfrak{g}_j(b) \right], P^\mathfrak{g}_j(x)\right], \theta P^\mathfrak{g}_j(y) \right], P^\mathfrak{g}_j(z) \right] \\
         &- \left[\left[ P^\mathfrak{g}_j(x), \left[\left[\theta P^\mathfrak{g}_j(b) , P^\mathfrak{g}_j(a)\right], \theta P^\mathfrak{g}_j(y) \right]\right], P^\mathfrak{g}_j(z) \right] \\
         &+ \left[\left[ P^\mathfrak{g}_j(x), \theta P^\mathfrak{g}_j(y) \right], \left[\left[ P^\mathfrak{g}_j(a), \theta P^\mathfrak{g}_j(b) \right], P^\mathfrak{g}_j(z) \right] \right]
    \end{align*}

    \item[(ii)] For any $a,x,y\in \mathfrak{g}$,
    \begin{align*}
        &\left[\left[ P^\mathfrak{g}_j(x), \theta P^\mathfrak{g}_j(y) \right], \left[\left[ P^\mathfrak{g}_j(x), \theta P^\mathfrak{g}_j(a) \right], P^\mathfrak{g}_j(x)\right]\right] \\
        &=\left[\left[ P^\mathfrak{g}_j(x), \left[\left[\theta P^\mathfrak{g}_j(y) , P^\mathfrak{g}_j(x)\right], \theta P^\mathfrak{g}_j(a) \right]\right], P^\mathfrak{g}_j(x) \right]
    \end{align*}
\end{itemize}
\end{lemma}
\begin{proof}
    (i) It follows from a straightforward application of Jacobi identity. (ii) The desired identity comes from a repeated application of assertion (i).
\end{proof}

Lemma \ref{p Lie triple identity} (i) will be referred to as the {\it TKK Jordan triple identity} since it coincides with the triple Jordan identity when applied in the setting of TKK Lie algebras associated with JB$^*$-triples.

\begin{lemma}\label{l perp reverse}
Let $(\mathfrak{g},\theta)$ be a TKK Lie algebra, and let $z,w\in \mathfrak{g}_{\pm 1}$ be such that $w$ satisfies \eqref{eq strict Lie tripotents} with $[[ P^\mathfrak{g}_j(w), \theta P^\mathfrak{g}_j(w) ], P^\mathfrak{g}_j(z) ]=0,$ for $j=-1,1$. Then $w\perp z$.
\end{lemma}
\begin{proof}
    By the linearity of $\theta$, a straightforward application of the Jacobi identity guarantees that, for $j=-1,1$, 
    \begin{equation}\label{eq l perp reverse 1}
        [[ P^\mathfrak{g}_j(z), \theta P^\mathfrak{g}_j(w) ], P^\mathfrak{g}_j(w) ]=0.
    \end{equation}

    On the other hand, since $w$ satisfies \eqref{eq strict Lie tripotents}, we apply Lemma \ref{p Lie triple identity} (ii) and \eqref{eq l perp reverse 1} to deduce that
     \begin{align}\label{eq l perp reverse 2}
         [[ P^\mathfrak{g}_j(w), \theta P^\mathfrak{g}_j(z) ], P^\mathfrak{g}_j(w) ]&= [[ P^\mathfrak{g}_j(w), \theta P^\mathfrak{g}_j(z) ], [[ P^\mathfrak{g}_j(w), \theta P^\mathfrak{g}_j(w) ], P^\mathfrak{g}_j(w) ] ] \\\nonumber
         &= [[ P^\mathfrak{g}_j(w), \theta [[ P^\mathfrak{g}_j(z), \theta P^\mathfrak{g}_j(w) ], P^\mathfrak{g}_j(w) ] ], P^\mathfrak{g}_j(w) ]\\\nonumber&=0.
     \end{align}

     Therefore, combining \eqref{eq l perp reverse 1} and \eqref{eq l perp reverse 2}, we have that
    \begin{align*}\label{eq l perp reverse 3}
         [ P^\mathfrak{g}_j(w), \theta P^\mathfrak{g}_j(z) ]&=[ [[ P^\mathfrak{g}_j(w), \theta P^\mathfrak{g}_j(w) ], P^\mathfrak{g}_j(w) ] ], \theta P^\mathfrak{g}_j(z) ] \\&=
         -[[ P^\mathfrak{g}_j(w), \theta P^\mathfrak{g}_j(z) ],[ P^\mathfrak{g}_j(w), \theta P^\mathfrak{g}_j(w) ]]\\&
         -[[ \theta P^\mathfrak{g}_j(z), [P^\mathfrak{g}_j(w) , \theta P^\mathfrak{g}_j(w)]], P^\mathfrak{g}_j(w) ]\\
         &= [[ P^\mathfrak{g}_j(w), \theta P^\mathfrak{g}_j(w) ],[ P^\mathfrak{g}_j(w), \theta P^\mathfrak{g}_j(z) ]]\\
         &= -[[ \theta P^\mathfrak{g}_j(w), [P^\mathfrak{g}_j(w) , \theta P^\mathfrak{g}_j(z)]], P^\mathfrak{g}_j(w) ]\\
         & -[[[ P^\mathfrak{g}_j(w), \theta P^\mathfrak{g}_j(z)] ,  P^\mathfrak{g}_j(w)], \theta P^\mathfrak{g}_j(w) ] = 0
     \end{align*}
     We have just shown that $z$ and $w$ are orthogonal.
\end{proof}

As a consequence of the result above, we have that any tripotent in a TKK Lie algebra satisfying \eqref{eq strict Lie tripotents} is orthogonal to its image through the TKK involution.

\begin{lemma}\label{l w perp theta w}
    Let $(\mathfrak{g},\theta)$ be a TKK Lie algebra, and let $w\in \mathcal{U}_s(\mathfrak{g_{\pm1}})$. Then, $w\perp \theta w$.
\end{lemma}
\begin{proof}
    Since $w$ is a tripotent in $\mathfrak{g_\pm1}$ satisfying \eqref{eq strict Lie tripotents}, we have that
    \begin{align*}
        w &= [[w,\theta w], w] = \sum_{j=-1,1}[[P^\mathfrak{g}_{j}(w), \theta P^\mathfrak{g}_{j}(w)], P^\mathfrak{g}_{j}(w)] \\&+  \sum_{j=-1,1}[[P^\mathfrak{g}_{j}(w), \theta P^\mathfrak{g}_{j}(w)], P^\mathfrak{g}_{-j}(w)].
    \end{align*}

    Therefore, for $j=-1,1$,
         \begin{align*}
       P^\mathfrak{g}_{j}(w)&= P^\mathfrak{g}_{j}(w)+ [[P^\mathfrak{g}_{-j}(w), \theta P^\mathfrak{g}_{-j}(w)], P^\mathfrak{g}_{j}(w)].
    \end{align*}
    And hence, for any $j=-1,1$, $$[[P^\mathfrak{g}_{-j}(w), \theta P^\mathfrak{g}_{-j}(w)], P^\mathfrak{g}_{j}(w)]=[[\theta P^\mathfrak{g}_{j}(\theta w), P^\mathfrak{g}_{j}(\theta w)], P^\mathfrak{g}_{j}(w)]=0.$$ An application of  Lemma \ref{l perp reverse} concludes the proof.
\end{proof}

The next lemma collects some facts that will be useful later.

\begin{lemma}\label{l trick perp}
Let $(\mathfrak{g},\theta)$ be a TKK Lie algebra,  let $x,y,w\in \mathfrak{g}_{\pm 1}$ be such that $y-x\perp x$, and let $j=-1, 1$. The following assertions hold.
\begin{itemize}
    \item[(i)] $[P^\mathfrak{g}_{j}(y), \theta P^\mathfrak{g}_{j}(x)] = [P^\mathfrak{g}_{j}(x), \theta P^\mathfrak{g}_{j}(x)].$
   \item[(ii)] 
        $[P^\mathfrak{g}_{j}(w-x), \theta P^\mathfrak{g}_{j}(x)]=[P^\mathfrak{g}_{j}(w-y), \theta P^\mathfrak{g}_{j}(x)]$.
    \item[(iii)] If $x$ and $y$ satisfy \eqref{eq strict Lie tripotents}, then
        $$P^\mathfrak{g}_{j}(x)= \left[ [P^\mathfrak{g}_{j}(y), \theta P^\mathfrak{g}_{j}(x)],P^\mathfrak{g}_{j}(y)\right].$$  
\end{itemize}
    
\end{lemma}

\begin{proof}
    (ii) By applying  (i) above, it can be deduced that, for $j=-1, 1$,
    \begin{align*}
        [P^\mathfrak{g}_{j}(w-x), \theta P^\mathfrak{g}_{j}(x)]=&[P^\mathfrak{g}_{j}(w), \theta P^\mathfrak{g}_{j}(x)]- [P^\mathfrak{g}_{j}(x), \theta P^\mathfrak{g}_{j}(x)]\\
        =& [P^\mathfrak{g}_{j}(w), \theta P^\mathfrak{g}_{j}(x)] -[P^\mathfrak{g}_{j}(y), \theta P^\mathfrak{g}_{j}(x)]\\
        =&[P^\mathfrak{g}_{j}(w-y), \theta P^\mathfrak{g}_{j}(x)].
    \end{align*}

    (iii) By assertion (i) and the Jacobi identity,  we have that, for $j=-1, 1$,
    \begin{align*}
        P^\mathfrak{g}_{j}(x)=& \left[ [P^\mathfrak{g}_{j}(x), \theta P^\mathfrak{g}_{j}(x)],P^\mathfrak{g}_{j}(x)\right]= \left[ [P^\mathfrak{g}_{j}(y), \theta P^\mathfrak{g}_{j}(x)],P^\mathfrak{g}_{j}(x)\right]\\
        =& -\left[ [\theta P^\mathfrak{g}_{j}(x), P^\mathfrak{g}_{j}(x)],P^\mathfrak{g}_{j}(y)\right] - \left[ [P^\mathfrak{g}_{j}(x), P^\mathfrak{g}_{j}(y)], \theta P^\mathfrak{g}_{j}(x)\right]\\
        =& \left[ [P^\mathfrak{g}_{j}(x), \theta P^\mathfrak{g}_{j}(x)],P^\mathfrak{g}_{j}(y)\right]
        =\left[ [P^\mathfrak{g}_{j}(y), \theta P^\mathfrak{g}_{j}(x)],P^\mathfrak{g}_{j}(y)\right].
         \end{align*}
    \end{proof}

 The concept of orthogonality amongst tripotents will allow for establishing a partial order relation on $\mathcal{U}_s(\mathfrak{g}_{\pm 1})$, as shown in Theorem \ref{t_partialorder}.

\begin{definition}\label{eq Lie partial order}
    Let $z,w\in \mathcal{U}_s(\mathfrak{g}_{\pm 1})$ be tripotents in a TKK Lie algebra $\mathfrak{g}=\mathfrak{g}_{-1}\oplus\mathfrak{g}_0\oplus \mathfrak{g}_1$ with involution $\theta$. We say that $z$ is {\rm less than or equal to}  $w$, written  $z\leq w$, if $w-z\in\mathcal{U}_s(g_{\pm 1})$ and $w-z\perp z$. 
\end{definition}

\begin{lemma}\label{l order in extremes}
Let $\mathfrak{g}=\mathfrak{g}_{-1}\oplus\mathfrak{g}_0\oplus \mathfrak{g}_1$ be a TKK Lie algebra with involution $\theta$. Let $x,y\in\mathcal{U}_s(\mathfrak{g}_{\pm 1})$ be such that $x\leq y$. Then, $P^\mathfrak{g}_j(x)\leq P^\mathfrak{g}_j(y)$, for  $j=-1,1$.
\end{lemma}

\begin{proof}
It is clear that $P^\mathfrak{g}_j(y)-P^\mathfrak{g}_j(x)=P^\mathfrak{g}_j(y-x)\in \mathfrak{g}_{\pm1}$, for $j=-1,1$. On the other hand, since $y-x\in \mathcal{U}_s(\mathfrak{g}_{\pm1})$, for $j=-1,1$, we have
$$
\left[ \left[  P^\mathfrak{g}_j(y-x), \theta P^\mathfrak{g}_j(y-x)  \right] , P^\mathfrak{g}_j(y-x) \right]  = P^\mathfrak{g}_j(y-x).
$$

Moreover, for $k=-1,1$ we have that
    \begin{align*}
        \left[ \left[  P^\mathfrak{g}_k(P^\mathfrak{g}_j(y-x)), \theta P^\mathfrak{g}_k(P^\mathfrak{g}_j(y-x))  \right] , P^\mathfrak{g}_k(P^\mathfrak{g}_j(y-x)) \right]  = 0 = P^\mathfrak{g}_k(P^\mathfrak{g}_j(y-x)),
    \end{align*}
    whenever $k\neq j$, and 
    \begin{align*}
        &\left[ \left[  P^\mathfrak{g}_k(P^\mathfrak{g}_j(y-x)), \theta P^\mathfrak{g}_k(P^\mathfrak{g}_j(y-x))  \right] , P^\mathfrak{g}_k(P^\mathfrak{g}_j(y-x)) \right]  \\&= \left[ \left[  P^\mathfrak{g}_j(y-x), \theta P^\mathfrak{g}_j(y-x)  \right] , P^\mathfrak{g}_j(y-x) \right]  = P^\mathfrak{g}_j(y-x) = P^\mathfrak{g}_k(P^\mathfrak{g}_j(y-x)),
    \end{align*}
    provided that $k=j$. We have just shown that $P^\mathfrak{g}_j(y-x)\in \mathcal{U}_s(\mathfrak{g}_{\pm1})$ for $j=-1,1$.

    Finally, we claim that, for any $j=-1,1$, $P^\mathfrak{g}_j(y-x)$ and $P^\mathfrak{g}_j(x)$ are orthogonal elements. Indeed, for any $k=-1,1$, we have that
    $$
     \left[  P^\mathfrak{g}_k(P^\mathfrak{g}_j(y-x)), \theta P^\mathfrak{g}_k(P^\mathfrak{g}_j(x))  \right] =0,
    $$ if $k\neq j$, and if $k=j$
    $$
     \left[  P^\mathfrak{g}_k(P^\mathfrak{g}_j(y-x)), \theta P^\mathfrak{g}_k(P^\mathfrak{g}_j(x))  \right] =  \left[  P^\mathfrak{g}_j(y-x), \theta P^\mathfrak{g}_j(x)  \right]=0,
    $$
    since $y-x\perp x$. Therefore, $P^\mathfrak{g}_j(x)\leq P^\mathfrak{g}_j(y)$, for any $j=-1,1$.
    \end{proof}

\begin{theorem}\label{t_partialorder}
     Let $(\mathfrak{g},\theta)$ be a TKK Lie algebra. The relation $\leq$  in Definition {\rm \ref{eq Lie partial order}} is a partial order on the subset $\mathcal{U}_s(\mathfrak{g}_{\pm 1})$ of tripotents in $\mathfrak{g}$.
\end{theorem}
\begin{proof}
Reflexivity. Let $z\in \mathcal{U}_s(\mathfrak{g}_{\pm 1})$. It is clear that $z\leq z$, since $z-z=0\in\mathcal{U}_s(\mathfrak{g}_{\pm 1})$ is a tripotent orthogonal to every tripotent:
    $$[ P^\mathfrak{g}_j(z-z), \theta P^\mathfrak{g}_j(z) ]= [ P^\mathfrak{g}_j(0), \theta P^\mathfrak{g}_j(z) ]=0.$$
    
Anti-symmetry. Let $z,w\in \mathcal{U}_s(\mathfrak{g}_{\pm 1})$, and suppose $z\leq w$ and $w\leq z$. We shall prove that $w=z$. On the one hand, it is clear that $P^\mathfrak{g}_0(w-z)=0$. 

On the other hand, since $w-z\in \mathcal{U}_s(\mathfrak{g}_{\pm 1})$, it follows that for each $j=-1, 1$, we have
    \begin{align*}
        P^\mathfrak{g}_j(w-z)=& \left[ [P^\mathfrak{g}_{j}(w-z), \theta P^\mathfrak{g}_{j}(w-z)],P^\mathfrak{g}_{j}(w-z)\right]\\
        =& \left[ [P^\mathfrak{g}_{j}(w-z), \theta P^\mathfrak{g}_{j}(w)]-  [P^\mathfrak{g}_{j}(w-z), \theta P^\mathfrak{g}_{j}(z)],P^\mathfrak{g}_{j}(w-z)\right]\\
        =& \left[ [P^\mathfrak{g}_{j}(w-z), \theta P^\mathfrak{g}_{j}(w)],P^\mathfrak{g}_{j}(w-z)\right]\\
        =& \left[ -[P^\mathfrak{g}_{j}(z-w), \theta P^\mathfrak{g}_{j}(w)],P^\mathfrak{g}_{j}(w-z)\right]=0.
    \end{align*}
    Therefore, $w-z= P^\mathfrak{g}_{-1}(w-z)+P^\mathfrak{g}_{0}(w-z)+P^\mathfrak{g}_{1}(w-z)=0$, as desired.\smallskip

Transitivity.  Let $x,y,w\in\mathcal{U}_s(\mathfrak{g}_{\pm 1})$ be such that $x\leq y$ and $y\leq w$. We shall prove that $x\leq w$. 

Firstly we show that $w-x$ and $x$ are orthogonal elements in $\mathfrak{g}_{\pm1}$. Indeed, it follows from Lemma \ref{l trick perp} (ii) and (iii)  that, for  $j=-1,1$,
    \begin{align*}
    \left[ P^\mathfrak{g}_{j}(w-x), \theta P^\mathfrak{g}_{j}(x)\right] =& \left[ P^\mathfrak{g}_{j}(w-y), \theta P^\mathfrak{g}_{j}(x)\right]\\
    = & \left[ P^\mathfrak{g}_{j}(w-y), \left[ [\theta P^\mathfrak{g}_{j}(y), P^\mathfrak{g}_{j}(x)],\theta P^\mathfrak{g}_{j}(y)\right]\right]\\
    =& \left[ [\theta P^\mathfrak{g}_{j}(y), P^\mathfrak{g}_{j}(w-y)],[\theta P^\mathfrak{g}_{j}(y), P^\mathfrak{g}_{j}(x)]\right] \\ & +
     \left[ [P^\mathfrak{g}_{j}(w-y), [\theta P^\mathfrak{g}_{j}(y), P^\mathfrak{g}_{j}(x) ]],\theta P^\mathfrak{g}_{j}(y)\right]\\
      =&\left[ [P^\mathfrak{g}_{j}(w-y), [\theta P^\mathfrak{g}_{j}(y), P^\mathfrak{g}_{j}(x) ]],\theta P^\mathfrak{g}_{j}(y)\right]\\
      =&\left[ 0, \theta P^\mathfrak{g}_{j}(y)\right] = 0.
    \end{align*}
We have just proved that $w-x\perp x$.\smallskip

We now claim that $w-x$ satisfies identities \eqref{eq strict Lie tripotents}. For this to be shown, we first note that, for $j=-1,1$,

    \begin{align*}
        \left[ P^\mathfrak{g}_{j}(w-x), \theta P^\mathfrak{g}_{j}(w-x)\right]=& \left[ P^\mathfrak{g}_{j}(w-x), \theta P^\mathfrak{g}_{j}(w)\right] - \left[ P^\mathfrak{g}_{j}(w-x), \theta P^\mathfrak{g}_{j}(x)\right]\\
        =& \left[ P^\mathfrak{g}_{j}(w-x), \theta P^\mathfrak{g}_{j}(w)\right]\\
        =& \left[ P^\mathfrak{g}_{j}(w), \theta P^\mathfrak{g}_{j}(w)\right]- \left[ P^\mathfrak{g}_{j}(x), \theta P^\mathfrak{g}_{j}(w)\right].
         \end{align*}
         Hence, by Lemma \ref{l trick perp} (i), for  $j=-1,1$,
     $$\left[ P^\mathfrak{g}_{j}(w-x), \theta P^\mathfrak{g}_{j}(w-x)\right] =\left[ P^\mathfrak{g}_{j}(w), \theta P^\mathfrak{g}_{j}(w)\right]- \left[ P^\mathfrak{g}_{j}(x), \theta P^\mathfrak{g}_{j}(x)\right].$$

This equality, together with Lemma \ref{l perp} (i), yields, for  $j=-1,1$,
    \begin{align*}
        [ [P^\mathfrak{g}_{j}(w-x), \theta P^\mathfrak{g}_{j}(w-x)]&,P^\mathfrak{g}_{j}(w-x)]= \left[ [P^\mathfrak{g}_{j}(w), \theta P^\mathfrak{g}_{j}(w)], P^\mathfrak{g}_{j}(w-x)\right]\\ &- \left[ [P^\mathfrak{g}_{j}(x), \theta P^\mathfrak{g}_{j}(x)], P^\mathfrak{g}_{j}(w-x)\right]\\
        &=\left[ [P^\mathfrak{g}_{j}(w), \theta P^\mathfrak{g}_{j}(w)], P^\mathfrak{g}_{j}(w-x)\right].
        \end{align*}
Therefore, by Lemma \ref{l trick perp} (i), 
        \begin{align*}
        [ [P^\mathfrak{g}_{j}(w-x), \theta P^\mathfrak{g}_{j}(w-x)],P^\mathfrak{g}_{j}(w-x)]
        =& P^\mathfrak{g}_{j}(w) - \left[ [P^\mathfrak{g}_{j}(w), \theta P^\mathfrak{g}_{j}(w)], P^\mathfrak{g}_{j}(x)\right]\\
        =& P^\mathfrak{g}_{j}(w) + \left[ [\theta P^\mathfrak{g}_{j}(w), P^\mathfrak{g}_{j}(x)], P^\mathfrak{g}_{j}(w)\right]\\
        =&  P^\mathfrak{g}_{j}(w) + \left[ [\theta P^\mathfrak{g}_{j}(x), P^\mathfrak{g}_{j}(x)], P^\mathfrak{g}_{j}(w)\right]\\
        =& P^\mathfrak{g}_{j}(w) - \left[ [P^\mathfrak{g}_{j}(w), \theta P^\mathfrak{g}_{j}(x)], P^\mathfrak{g}_{j}(x)\right]\\
        =& P^\mathfrak{g}_{j}(w) - \left[ [P^\mathfrak{g}_{j}(x), \theta P^\mathfrak{g}_{j}(x)], P^\mathfrak{g}_{j}(x)\right]\\
        =& P^\mathfrak{g}_{j}(w) - P^\mathfrak{g}_{j}(x)= P^\mathfrak{g}_{j}(w-x),
    \end{align*}
    for  $j=-1, 1$, as desired.\smallskip

    It is left to prove that $w-x$ is a tripotent. In order to do that, observe that
    \begin{align*}
        w-x&=[[w-x, \theta(w-x)], w-x] + [[x, \theta(w-x)], w]\\&+ [[x, \theta x], w-x] + [[w-x, \theta x], w] + [[w-x, \theta(w-x)], x].
    \end{align*}

    By orthogonality and Lie arithmetic, we have that
    \begin{align*}
        [[x, \theta(w-x)], w]&= \sum_{j=-1,1}[[P^\mathfrak{g}_{j}(x), \theta P^\mathfrak{g}_{j}(w-x)], w] \\&+  \sum_{j=-1,1}[[P^\mathfrak{g}_{j}(x), \theta P^\mathfrak{g}_{-j}(w-x)], w]=0.
    \end{align*}
    It can be similarly proved that $[[w-x, \theta x], w]=0$.\smallskip

    On the other hand, by Lemma \ref{l perp} (i),
    $$
    [[x,\theta x],w-x] = \sum_{j=-1,1}[[P^\mathfrak{g}_{j}(x), \theta P^\mathfrak{g}_{j}(x)], P^\mathfrak{g}_{-j}(w-x)].
    $$
    Since $x,w\in \mathcal{U}_s(\mathfrak{g_{\pm 1}})$, by Lemma \ref{l w perp theta w} and Lemma \ref{l trick perp} (i), we have that
    \begin{align*}
        [ [P^\mathfrak{g}_{j}(x), \theta P^\mathfrak{g}_{j}(x)],P^\mathfrak{g}_{-j}(w-x)] &= [ [P^\mathfrak{g}_{j}(x), \theta P^\mathfrak{g}_{j}(x)],\theta P^\mathfrak{g}_{j}(\theta(w-x))]
        \\&= [ [P^\mathfrak{g}_{j}(x), \theta P^\mathfrak{g}_{j}(x)],\theta P^\mathfrak{g}_{j}(\theta w)]\\&=
         [ [P^\mathfrak{g}_{j}(w), \theta P^\mathfrak{g}_{j}(x)],\theta P^\mathfrak{g}_{j}(\theta w)]=0,
    \end{align*}
    for $j=-1,1$. Therefore, $[[x,\theta x],w-x]=0$. Similar arguments give that, for $j=-1,1$, 
    $$[ [P^\mathfrak{g}_{j}(w-x), \theta P^\mathfrak{g}_{j}(w-x)],P^\mathfrak{g}_{-j}(x)]=0,
    $$
    and hence 
    $$[[w-x,\theta(w-x)],w-x]=0.
    $$

    Combining the previous assertions, we deduce that $$ w-x=[[w-x, \theta(w-x)], w-x],$$ as claimed. Consequently, $x\leq w$ in $\mathcal{U}_s(\mathfrak{g}_{\pm 1})$.

\end{proof}

The partial order on the set of tripotents of a JB$^*$-triple $V$ aligns with the Lie-algebraic partial order of Definition \ref{eq Lie partial order} on the subset $\mathcal{U}_s(\fL(V)_{\pm 1})$ of tripotents of its TKK Lie algebra $\fL(V)$. The following proposition illustrates that fact.

\begin{proposition}\label{p order characterisation}
Let $V$ be a JB*-triple, let $\mathfrak{L}(V)=\mathfrak{L}(V)_{-1}\oplus\mathfrak{L}(V)_0\oplus \mathfrak{L}(V)_1$ be the  TKK Lie algebra of $V$, and let 
$z,w\in \mathcal{U}(V)$ be tripotents in $V$. Then, $z\leq w$ in the JB*-triple $V$ if and only if $z\leq w$ in the TKK Lie algebra $\mathfrak{L}(V)$. 
\end{proposition}

\begin{proof}
    Let  $z,w$ be tripotents in $V$. By means of the usual identification of $V$ and $\mathfrak{L}(V)_{-1}$,  $z=P^{\mathfrak{L}(V)}_{-1}(z)$ and $w=P^{\mathfrak{L}(V)}_{-1}(w)$ are also considered as elements of $\mathfrak{L}(V)$. By Proposition \ref{p tripotent characterisation} and comments below, we have that $z$ and $w$ are tripotents in the subset $\mathcal{U}_s(\mathfrak{L}(V)_{\pm 1})$. 

    Suppose that $z\leq w$ in $V$, that is, $w-z$ is a tripotent in $V$ orthogonal to $z$. Hence, by Proposition \ref{p tripotent characterisation}, we have that $w-z\in \mathcal{U}_s(\mathfrak{L}(V)_{\pm 1})$. Moreover, it follows from Proposition \ref{p orth characterisation} that $w-z\perp z$ in $\mathfrak{L}(V)$.

    On the other hand, if $z$ and $w$ are regarded as tripotents in the poset $\mathcal{U}_s(\mathfrak{L}(V)_{\pm 1})$ and are such that $z\leq w$, then the tripotent $w-z$ is orthogonal to $z$ in $\mathfrak{L}(V)$, by definition. Hence, by Proposition \ref{p tripotent characterisation}, $w-z\in \mathcal{U}(V)$. Furthermore, by Proposition \ref{p orth characterisation}, we have that $w-z$ and $z$ are orthogonal in  $V$. We have just proved that $z\leq w$ in the JB*-triple $V$.
\end{proof}

 Let $(\mathfrak{g}, \theta)$ and $(\mathfrak{h}, \theta')$ be two TKK Lie algebras. Consider a bijection $\varphi: \mathcal{U}_s(\mathfrak{g}_{\pm1})\to \mathcal{U}_s(\mathfrak{h}_{\pm1})$. The mapping $\varphi$ is said to preserve the partial order $\leq$ on $\mathcal{U}_s(\mathfrak{g}_{\pm1})$ in one direction if, for any two elements $z,w\in \mathcal{U}_s(\mathfrak{g}_{\pm1})$ such that $z\leq w$, then $\varphi(z)\leq \varphi(w)$ in $\mathcal{U}_s(\mathfrak{h}_{\pm1})$. The mapping $\varphi$ preserves the partial order in both directions if the equivalence $z\leq w \iff \varphi(z)\leq \varphi(w)$ holds, for any $z,w\in \mathcal{U}_s(\mathfrak{g}_{\pm1})$. An {\it order isomorphism} between posets is a bijection preserving partial order in both directions. Since the TKK involution $\theta$ of a TKK Lie algebra $\mathfrak{g}$ preserves the elements of $\mathcal{U}_s(\mathfrak{g}_{\pm1})$ and the orthogonality amongst them, it is not difficult to see that $\theta$ is always an order isomorphism on the poset $\mathcal{U}_s(\mathfrak{g}_{\pm1})$. \smallskip

Lemma \ref{dsum} considers TKK Lie algebras associated with atomic JB*-triples. The following theorem is a Lie-algebraic analogue of \cite[Theorem 6.1]{FriPer} proved by Y. Friedmann and A.M. Peralta in the Jordan triple setting.

\begin{theorem}\label{t_atomic}
Let $V= \bigoplus_{\gamma\in\Gamma}^\infty  V_\gamma$ and $W= \bigoplus_{\lambda\in \Lambda}^\infty  W_\lambda$ be atomic JBW*-triples, where all summands are Cartan factors of rank greater than or equal to $2$. Let  $(\mathfrak{L}(V),\theta)$ and $(\mathfrak{L}(W),\theta')$ be the   TKK Lie algebras associated with $V$ and $W$, respectively. Let $\Delta: \mathcal{U}_s(\mathfrak{L}(V)_{\pm 1}) \to \mathcal{U}_s(\mathfrak{L}(W)_{\pm 1})$ be a graded order isomorphism commuting with involutions which preserves orthogonality  and is continuous at some tripotent $z=(z_\gamma)\in \mathcal{U}_s(\mathfrak{L}(V)_{\pm 1})$ such that one of the following conditions holds:
\begin{itemize}
\item [(i)]  $P_{-1}^{\mathfrak{L}(V)}(z_\gamma)\neq 0$, for all $\gamma\in \Gamma$; 
\item [(ii)]$P_{1}^{\mathfrak{L}(V)}(z_\gamma)\neq 0$, for all $\gamma\in \Gamma$. 
\end{itemize}
Then, there exists a real-linear graded isomorphism  $T:\mathfrak{L}(V)\to \mathfrak{L}(W)$ extending $\Delta$, that is, for all $z\in \mathcal{U}_s(\mathfrak{L}(V)_{\pm 1})$,  $T(z)=\Delta(z)$. 
\end{theorem}

\begin{proof}
Let $\mathfrak{L}(V)=\mathfrak{L}(V)_{-1}\oplus\mathfrak{L}(V)_{0}\oplus\mathfrak{L}(V)_{1}$ and $\mathfrak{L}(W)=\mathfrak{L}(W)_{-1}\oplus\mathfrak{L}(W)_{0}\oplus\mathfrak{L}(W)_{1}$ be the gradings of the TKK Lie algebras $\mathfrak{L}(V)$ and $\mathfrak{L}(W)$, respectively.

Notice that, by the definition of the sets $\mathcal{U}_s(\mathfrak{L}(V)_{\pm 1})$ and $ \mathcal{U}_s(\mathfrak{L}(W)_{\pm 1})$ of tripotents,  we have that, for $j=-1,1$, 
\begin{equation}\label{301}
\mathcal{U}_s(\mathfrak{L}(V)_{\pm 1})\cap \mathfrak{L}(V)_{j}=\mathcal{U}(\mathfrak{L}(V)_{j})
\end{equation}
and
\begin{equation}\label{302}
\mathcal{U}_s(\mathfrak{L}(W)_{\pm 1})\cap \mathfrak{L}(W)_{j}=\mathcal{U}(\mathfrak{L}(W)_{j})
\end{equation}
(cf. Definition \ref{d tripotent} and Proposition \ref{p tripotent characterisation}).
 
 Let $\Delta=\Delta_{-1}+\Delta_0+\Delta_1$ be the grading of $\Delta$, where $\Delta_0=0$  and $\Delta_j: \mathcal{U}_s(\mathfrak{L}(V)_{\pm 1})\cap \mathfrak{L}(V)_{j}\to \mathcal{U}_s(\mathfrak{L}(W)_{\pm 1})\cap \mathfrak{L}(W)_{j}$, for $j=-1,1$.. In particular, by Proposition \ref{p tripotent characterisation} and equalities \eqref{301}, \eqref{302}, the mapping  $\Delta_{-1}:\mathcal{U}(\mathfrak{L}(V)_{-1})\to \mathcal{U}(\mathfrak{L}(W)_{-1})$ can and will be seen as a bijection between the set of tripotents of $V=\mathfrak{L}(V)_{-1}$ and $W=\mathfrak{L}(W)_{-1}$. The mapping $\Delta_{-1}:\mathcal{U}(\mathfrak{L}(V)_{-1})\to \mathcal{U}(\mathfrak{L}(W)_{-1})$  is an orthogonality preserving order isomorphism.  Moreover,  assume that (i) is satisfied. Then there exists a $u=(u_\gamma)\in \mathcal{U}(\mathfrak{L}(V)_{-1})$ with $u_\gamma\neq 0$, for all $\gamma$, at which $\Delta_{-1}$ is continuous. In fact, consider $u$ as the identification of $P_{-1}(z)$ as an element in $V$. If we assume (ii) to be satisfied instead, the same argument applies to $\theta z$. The continuity of the restriction is guaranteed since $\Delta$ commutes with the involutions (which are continuous). 

It is now the case that $\Delta_{-1}$ satisfies \cite[Theorem 6.1]{FriPer} and, consequently, there exists a real-linear triple isomorphism $\phi:V\to W$ such that $\phi(x)=\Delta_{-1}(x)$, for all  tripotents $x\in \mathfrak{L}(V)_{-1}$. (See also Proposition \ref{p orth characterisation} above.)

Observe that, since $\Delta$ commutes with involutions, it follows that $\phi$ also extends $\Delta_1$. To be precise, the  mapping $\psi\colon \mathfrak{L}(V)_1\to\mathfrak{L}(W)_1$ defined , for all $y\in \mathfrak{L}(V)_1$, by 
$$\psi(y)= \psi(\theta x)=\theta \phi (x),
$$
where $x\in \mathfrak{L}(V)_{-1}$ is the unique element in $\mathfrak{L}(V)_{-1}$ such that $y=\theta x$, extends $\Delta_1$.
To avoid cumbersome notation, we will use $\phi$ to denote both extensions.

We define a graded mapping $T:\mathfrak{L}(V)\to \mathfrak{L}(W)$ such that 
    $$
    T(x,h,y)=(\phi(x), \tilde{\phi}(h), \phi(y)), \qquad (x,h,y)\in \mathfrak{L}(V),
    $$
where $\tilde{\phi}(h):=\phi\circ h\circ \phi^{-1}$, for each $h\in \mathfrak{L}(V)_0. $
 
We prove next that $T:\mathfrak{L}(V)\to \mathfrak{L}(W)$ is a $\mathbb{R}$-linear isomorphism extending $\Delta$.

Firstly, we show that $\tilde{\phi}:\mathfrak{L}(V)_0\to \mathfrak{L}(W)_0$  $\tilde{\phi}$ is a real-linear bijection. Indeed, since $\phi$ is real-linear, for each $h\in \mathfrak{L}(V)_0$ and each $\lambda\in \mathbb{R}$, $ \tilde{\phi}(\lambda h) = \phi\circ \lambda h \circ \phi^{-1} =  \lambda \phi\circ h \circ \phi^{-1}=\lambda \tilde{\phi}(h)$, yielding that $\tilde{\phi}$ is real-linear. 

On the other hand, given some $h'\in \mathfrak{L}(W)_0$, there exists $h\in \mathfrak{L}(V)_0$ such that $\tilde{\phi}(h)=h'$. Indeed, it is enough to set $h=\phi^{-1}\circ h'\circ \phi$. Hence, $\tilde{\phi}$ is surjective.

It now remains to show that $\tilde{\phi}$ is injective.  
Let $h=\sum_{j=1}^n a_j\bo b_j$ and $k=\sum_{l=1}^m c_l\bo d_l$  be arbitrary elements in $\mathfrak{L}(V)_0$, and suppose that $\tilde{\phi}(h)=\tilde{\phi}(k)$. Then, since $\phi$ preserves the triple product, we have that
\begin{align*}
      \sum_{j=1}^n \phi(a_j)\bo \phi(b_j)=&\phi \circ \sum_{j=1}^n a_j\bo b_j \circ \phi^{-1} = \phi\circ h \circ \phi^{-1}= \tilde{\phi}(h)\\=& \tilde{\phi}(k)=\phi\circ k\circ \phi^{-1}=\sum_{l=1}^m \phi(c_l)\bo \phi(d_l).
\end{align*}
It is also the case that, since $\phi$ is a bijection, we have, for any $w\in \mathfrak{L}(W)_{-1}$, that there exists a unique $z\in \mathfrak{L}(V)_{-1}$  such that $\phi(z)=w$. Hence, by the equality above,
\begin{align*}
        \sum_{j=1}^n \phi(a_j)\bo \phi(b_j)(w) =& \sum_{j=1}^n \phi(a_j)\bo \phi(b_j)(\phi(z))=\sum_{l=1}^m \phi(c_l)\bo \phi(d_l)(\phi(z)). 
    \end{align*}
In other words, for all $z\in \mathfrak{L}(V)_{-1},$
$$\phi (\sum_{j=1}^n\{a_j,b_j,z\}) = \phi (\sum_{l=1}^m\{c_l,d_l,z\}),$$
yielding finally that 
$$h=\sum_j a_j\bo b_j = \sum_l c_l\bo d_l=k,$$ as required.

  It now follows immediately from the grading that  $T$ is a real-linear bijection.

    We show now that $T$ preserves the Lie bracket.  Consider any two elements $(x,h,y),(u,k,v)\in \mathfrak{L}(V)$. 

    Claim 1. $\phi(h(u))=\tilde{\phi}(h)(\phi(u))$.\\
   This is straightforward since $\tilde{\phi}(h)(\phi(u))=\phi\circ h \circ \phi^{-1}(\phi(u))=\phi(h(u)).$
    
    Claim 2. $\tilde{\phi}$ is multiplicative, that is, $\tilde{\phi}(hk)=\tilde{\phi}(h)\tilde{\phi}(k)$.\\
    Indeed, 
    $$\tilde{\phi}(hk)=\phi\circ hk\circ \phi^{-1}=\phi\circ h\circ \phi^{-1}\circ \phi \circ k\circ \phi^{-1}=\tilde{\phi}(h)\tilde{\phi}(k),$$
which proves the claim.

    Claim 3. If $k=\sum_{l=1}^m c_l\bo d_l$, then  $\phi(k^\natural(y))=\tilde{\phi}(k)^\natural(\phi(y))$.\\
    We have
    \begin{align*}
        \phi(k^\natural(y))&=\phi(\sum_{l=1}^m d_l\bo c_l(y))=\sum_{l=1}^m \phi(d_l)\bo \phi(c_l)(\phi(y))\\&=(\phi\circ \sum_{l=1}^m c_l\bo d_l \circ \phi^{-1})^\natural(\phi(y))=\tilde{\phi}(k)^\natural(\phi(y)),
    \end{align*}
as desired.

    Combining the three claims above, we have
    \begin{align*}
        T[(x,y,h),(u,k,v)]&= T(h(u)-k(x), [h,k]+x\bo v-u\bo y, k^\natural(y)-h^\natural(v)),\\&= (\phi(h(u)-k(x)), \phi\circ ([h,k]+x\bo v-u\bo y)\circ \phi^{-1},\\ & \phi(k^\natural(y)-h^\natural(v)))\\&=
        (\tilde{\phi}(h)(\phi(u))-\tilde{\phi}(k)(\phi(x)),\\& \tilde{\phi}(h)\tilde{\phi}(k)-\tilde{\phi}(k)\tilde{\phi}(h)+\phi\circ x\bo v\circ \phi^{-1} - \phi\circ u\bo y\circ \phi^{-1},\\& \tilde{\phi}(k)^\natural(\phi(y))-\tilde{\phi}(h)^\natural(\phi(v)) )\\&=
        (\tilde{\phi}(h)(\phi(u))-\tilde{\phi}(k)(\phi(x)),\\& [\tilde{\phi}(h),\tilde{\phi}(k)]+\phi(x)\bo \phi(v) - \phi(u)\bo \phi(y),\\& \tilde{\phi}(k)^\natural(\phi(y))-\tilde{\phi}(h)^\natural(\phi(v)) )\\
        &=[T(x,h,y), T(u,k,v)],
    \end{align*}
    yielding that $T[(x,y,h),(u,k,v)]=[T(x,y,h),T(u,k,v)]$, which shows that $T$ preserves the Lie bracket.

    It only remains to prove that $T$ commutes with involutions. For $(x,h,y)\in \mathfrak{L}(V)$,
    \begin{align*}
        T\theta(x,h,y)=T(y,-h^\natural,x)=(\phi(y), \phi\circ (-h^\natural)\circ \phi^{-1}, \phi(x)).
    \end{align*}
    Observe that 
    \begin{align*}
        \phi\circ (-h^\natural)\circ \phi^{-1}&=\phi\circ (-(\sum_j a_j\bo b_j)^\natural)\circ \phi^{-1}\\&=\phi\circ (-\sum_j b_j\bo a_j)\circ \phi^{-1}\\&=-\phi\circ (\sum_j b_j\bo a_j)\circ \phi^{-1}=\\&=-(\phi\circ h\circ \phi^{-1})^\natural.
    \end{align*}
Therefore, it is clear that $T\theta(x,h,y)=\theta'T(x,h,y)$.

It is easily seen that $T(z)=\Delta(z)$, for all tripotents $z$ in $\mathcal{U}_s(\mathfrak{L}(V)_{\pm 1})$, which concludes the proof.
\end{proof}

An immediate consequence of this theorem is the following corollary.

\begin{corollary}\label{c_atomic}
Let $V= \bigoplus_{\gamma\in\Gamma}^\infty  V_\gamma$ and $W= \bigoplus_{\lambda\in\Lambda}^\infty  W_\lambda$ be  atomic JBW*-triples, where all summands are Cartan factors of rank greater than or equal to $2$. Let  $(\mathfrak{L}(V),\theta)$ and $(\mathfrak{L}(W),\theta')$ be the   TKK Lie algebras associated with $V$ and $W$, respectively. Let $\Delta: \mathcal{U}_s(\mathfrak{L}(V)_{\pm 1}) \to \mathcal{U}_s(\mathfrak{L}(W)_{\pm 1})$ be a  negatively graded   order isomorphism commuting with involutions which  preserves orthogonality,  and is continuous at some tripotent $z=(z_\gamma)\in \mathcal{U}_s(\mathfrak{L}(V)_{\pm 1})$ such that one of the following conditions holds:
\begin{itemize}
\item [(i)]  $P_{-1}^{\mathfrak{L}(V)}(z_\gamma)\neq 0$, for all $\gamma\in \Gamma$; 
\item [(ii)]$P_{1}^{\mathfrak{L}(V)}(z_\gamma)\neq 0$, for all $\gamma\in \Gamma$. 
\end{itemize}
Then, there exists a real-linear  negatively graded isomorphism  $T:\mathfrak{L}(V)\to \mathfrak{L}(W)$ extending $\Delta$, that is, for all $z\in \mathcal{U}_s(\mathfrak{L}(V)_{\pm 1})$,  $T(z)=\Delta(z)$ . 
\end{corollary}

\begin{proof}
Let $\Delta: \mathcal{U}_s(\mathfrak{L}(V)_{\pm 1}) \to \mathcal{U}_s(\mathfrak{L}(W)_{\pm 1})$ be  negatively graded. Consider the mapping $\Delta': \mathcal{U}_s(\mathfrak{L}(V)_{\pm 1}) \to \mathcal{U}_s(\mathfrak{L}(W)_{\pm 1})$ defined by 
$\Delta':=\theta' \Delta$. 
Since the TKK involution $\theta'$ is a negatively graded continuous bijection preserving tripotents, orthogonality, and the partial ordering on  $\mathcal{U}_s(\mathfrak{L}(W)_{\pm1})$ in both directions, it follows that $\Delta'$ satisfies the conditions of Theorem \ref{t_atomic} (see Proposition \ref{p tripotent characterisation} and Remark \ref{remx2}). Therefore, by Theorem \ref{t_atomic}, there exists a real-linear graded isomorphism  $T':\mathfrak{L}(V)\to \mathfrak{L}(W)$ extending $\Delta'$. Consequently, the real-linear negatively graded isomorphism given by 
 $T:=\theta' T'$ extends $\Delta.$

\end{proof}

\quad\\

\textbf{Acknowledgements.} The authors wish to thank Professor Cho-Ho Chu and Dr. Pedro Lopes for valuable discussions during the preparation of the manuscript.

\qquad\\


\begin{thebibliography}{0}
\bibitem{Adamek}J. Adamek, H. Herrlich, G.E. Strecker,  Abstract and Concrete Categories: The Joy of Cats,  Dover Publications (2009).
\bibitem{batt} M. Battaglia, Order theoretic type decomposition of JBW*-triples, {\it Quart. J. Math.} {\bf 42} (1991), 129--147.
\bibitem{BarTi86} T. Barton, R.M. Timoney, Weak*-continuity of Jordan triple products and applications, {\it Math. Scand.} \textbf{59} (1986), 177--191.

\bibitem{CuetoBecFerPer} J. Becerra Guerrero, M. Cueto-Avellaneda, F.J. Fern\'andez-Polo, A.M. Peralta, On the extension of isometries between the unit spheres of a JBW$^*$-triple and a Banach space, {\it Journal of the Institute of Mathematics of Jussieu} \textbf{20} 1, 277-303 (2021)

\bibitem{BraKaUp78} R.B. Braun, W. Kaup, H. Upmeier, A holomorphic characterization of Jordan-C*-algebras, {\it Math. Z.} \textbf{161}, 277-290 (1978).


\bibitem{CavSmir}
D. M. Caveny and O. N. Smirnov,
Categories of Jordan structures and graded Lie algebras, {\it Comm. Algebra} \textbf{42} 186-202 (2014).

\bibitem{ChuBook} C-H. Chu, Jordan structures in geometry and analysis, Cambridge Univ. Press, Cambridge, 2012.


\bibitem{ChuOliveira} C-H. Chu, L. Oliveira, Tits-Kantor-Koecher Lie algebras of JB*-triples {\it Journal of Algebra} \textbf{512} 465-492 (2018).


\bibitem{Dang92} T. Dang, Real isometries between JB*-triples, {\it Proc. Amer. Math. Soc.} \textbf{114}, 971-980 (1992).

\bibitem{Di86} S. Dineen, The second dual of a JB*-triple system, In: Complex analysis, functional analysis and approximation theory (ed. by J. M\'ugica), 67-69, (North-Holland Math. Stud. 125), North-Holland, Amsterdam-New York, (1986).

\bibitem{EdFerHosPe2010} C.M. Edwards, F.J. Fern\'andez-Polo, C.S. Hoskin, A.M. Peralta, On the facial structure of the unit ball in a JB$^*$-triple, {\it J. Reine Angew. Math.} \textbf{641}, 123-144 (2010).


\bibitem{EdRu96} C.M. Edwards, G.T. R\"uttimann, Compact tripotents in bi-dual JB$^*$-triples, Math. Proc. Camb. Phil. Soc. \textbf{120}, 155-173 (1996).


\bibitem{ERFS} C.M. Edwards and G.T. R\"uttimann, On the facial structure of the unit balls in a JBW*-triple and its predual. {\it J. London Math. Soc.} {\bf 38} (1988), 317--332.


\bibitem{FriPer} Y.~{Friedman} and A.M.~{Peralta}. \newblock {Representation of symmetry transformations on the sets of tripotents of spin and Cartan factors.} \newblock {\it {Anal. Math. Phys.}}, \textbf{12} 37, (2022).

\bibitem{FriRu85} Y.~{Friedman} and B.~{Russo}. \newblock {Structure of the predual of a JBW$\sp*$-triple.} \newblock {\it {J. Reine Angew. Math.}}, \textbf{356} 67-89, 1985.

\bibitem{Gar} E. Garcia,
Tits-Kantor-Koecher algebras of strongly prime Hermitian Jordan pairs, {\it J. Algebra} \textbf{277} 559-571 (2004).

\bibitem{Helgason} S. Helgason, Differential geometry, Lie groups and symmetric spaces, Academic Press, San Diego, 1978.

\bibitem{Horn87} G. Horn, Characterization of the predual and ideal structure of a JBW*-triple, {\it Math. Scand.} \textbf{61}, no. 1, 117-133 (1987).


\bibitem{IsKaRo95} J.M. Isidro, W. Kaup, A. Rodr{\'\i}guez, On real forms of JB*-triples, {\it Manuscripta Math.} \textbf{86}, 311-335 (1995).

\bibitem{Jacob} N. Jacobson, {\it Structure and Representation of Jordan algebras}, Amer. Math. Soc. Providence, R.I. 1969.

\bibitem{k1} I.L. Kantor, {\it Classification of irreducible
transitive differential groups}, Dokl. Akad. Nauk SSSR
158 (1964) 1271-1274.

\bibitem{k2} I.L. Kantor, {\it Transitive differential groups and
invariant connections on homogeneous spaces}, Trudy Sem. Vecktor.
Tenzor. anal. 13 (1966) 310-398.


\bibitem{Ka83} W. Kaup, A Riemann Mapping Theorem for bounded symmentric domains in complex Banach spaces, {\it Math. Z.} \textbf{183}, 503-529 (1983).

\bibitem{KaUp77} W. Kaup, H. Upmeier, Jordan algebras and symmetric Siegel domains in Banach spaces, {\it Math. Z.} \textbf{157}, 179-200 (1977).

\bibitem{Koecher} M. Koecher, \textit{Imbedding of Jordan algebras into Lie algebras I}, Amer. J. Math. 89, Nr. 3, (1967) 787-816.


\bibitem{loos1} O. Loos, Jordan triple systems, R-spaces, and
bounded symmetric domains, {\it Bull. Amer. Math. Soc.} \textbf{77} 558-561 (1971).

\bibitem{MacLane} S. Mac Lane, {\it Categories for the Working Mathematician}, Graduate texts in Mathematics  \textbf{5} Springer-Verlag (1971).

\bibitem{my} K. Meyberg, {\it Jordan-Tripelsysteme und die
Koecher-Konstruktion von Lie-Algebren}, Math. Z. 115
(1970) 58-78.

\bibitem{Neher} E. Neher, 3-graded Lie algebras and Jordan pairs, {\it Non-associative algebra and its applications} (Oviedo, 1993) 296-299, Kluwer Acad.Publ. Dordrecht, 1994.


\bibitem{Riehl} E. Riehl,  \href{https://math.jhu.edu/~eriehl/context.pdf}{\it Category Theory in Context} (https://math.jhu.edu/~eriehl/context.pdf).


\bibitem{Sidd2007} A.A. Siddiqui, Average of two extreme points in JBW$^*$-triples, {\it Proc. Japan Acad.} Ser. A \textbf{83}, No. 9-10, 176-178 (2007).

\bibitem{tits} J. Tits, {\it Alg\`ebres alternatives, alg\`ebres de Jordan et alg\`ebres de Lie exceptionnelles I: Construction}, Nederl. Indag. Math. 28 (1966) 223-237.



\end{thebibliography}
\end{document}